# Development of an Efficient Formulation for Volterra's Equations of Motion for Multibody Dynamical Systems

Mohammad Hussein Yoosefian Nooshabadi, Hossein Nejat Pishkenari*

*Center of Excellence in Design, Robotics and Automation (CEDRA), Department of Mechanical Engineering, Sharif University of Technology, Tehran, Iran*

**Abstract**

In this paper, we present an efficient form of Volterra's equations of motion for both unconstrained and constrained multibody dynamical systems that include ignorable coordinates. The proposed method is applicable for systems with both holonomic and nonholonomic constraints. Firstly, based on the definition of ignorable coordinates, one of the motion constants (the generalized momentum vector corresponding to the ignorable coordinates) is dealt with as a constraint, which will be referred to as dynamical constraints. These constraints, along with ordinary constraints, namely kinematical constraints, are then used in the proposed method to derive motion equations. This approach gives the minimum number of equations needed to study the behavior of a dynamical system. Three simulation examples are provided to evaluate the proposed method and to compare it to existing methods. The first case study is a constrained dynamical system, which moves in two-dimensional space. The second one is an unconstrained multibody system including three connected rigid bodies. Finally, the last case study includes a cubic satellite that uses a deployable boom to move a mass to a desired location. The results of the numerical simulations are compared to the conventional methods and the better performance of the proposed method is demonstrated.

*Keywords:* Volterra's Equations of motion; Constrained mechanical system; Ignorable coordinates; Quasi-velocity; Dynamical constraints

## 1. Introduction

Deriving the equations of motion governing dynamical systems is the foremost step in the analysis of their behavior, and therefore their control. First attempts in deriving equations of motion via analytical approach trace back to 1788, when J. L. Lagrange introduced the notion of generalized coordinates. Lagrange categorized forces acting on a constrained dynamical system into two groups, that is, forces imposed by constraints and other external forces. By introducing generalized coordinates and using Lagrange multipliers, he tried to derive equations of motion of a system [1, 2]. Since the advent of the Lagrange equations, a variety of approaches emerged, all of which have their pros and cons [3]. Following Lagrange, Hamilton developed his well-known principle that could be used to derive equations of both discrete and continuous systems [4]. Based on the Hamiltonian, he tried to transform the second-order dynamical equations to first-order ones. In the1840's, Jacobi developed the well-known Jacobi integral, which generalizes the energy of a mechanical system [5, 6], based on which he proposed a new method of describing the motion of a system [7]. In the late 19$^{th}$, Routh developed an approach to eliminate the ignorable coordinates of an unconstrained mechanical system. Using the notion of generalized momentum and introducing the so-called Routhian function, he, similar to Hamilton, obtained the first-order form of the equations of motion [8]. J. W. Gibbs introduced the concept of quasi-velocities and used them to explore systems' dynamics [9]. In the Gibbs method, a scalar function is used, which is based on the accelerations of the bodies in the system and could be interpreted as acceleration energy (an equivalent term to kinetic energy in the Lagrange method) [9]. Based on the fundamental equations of motion, in 1898, Volterra derived a special form of the equations of motion in which quasi-velocities were employed but, unlike the Gibbs approach, accelerations were not required to be calculated, therefore less computational effort was needed [10]. G. A. Maggi proposed a Lagrange-based method for deriving ordinary differential equations (ODEs) governing a constrained system. Unlike the Lagrange method, these equations did not include Lagrange multipliers and therefore were simpler [11]. In the early 20$^{th}$, Hamel considered mechanical systems under nonlinear constraints [12]. In the late 20$^{th}$, T.R. Kane developed a new approach which, similar to the Gibbs approach, was based on quasi-velocities [13-15]. Finally, the most recent approach was developed by Udwadia and Kalaba in 2007, in which no auxiliary variables were employed [16-18].

*Corresponding Author. Room 253, Department of Mechanical Engineering, Sharif University of Technology, Azadi Avenue, Tehran, Iran. Tel: 0098-2166165543, Fax: 0098-2166000021, Email address: nejat@sharif.edu

Since virtually all dynamical systems could be analyzed using Lagrange equations, other methods should relatively be advantageous in some aspects to warrant interest. One prominent aspect of an approach is its performance during the numerical solution. When it comes to the numerical solution of the equations of constrained dynamical systems, two factors play key roles in assessing the performance of an approach, viz. numerical integration error, and the amount of time that the simulation takes to run. The numerical error occurs and accumulates because of the imprecise numerical integrating methods. Many researchers have tried to reduce the computational error [19-23]. The amount of numerical error pertains to two major factors. The first one is the order of equations of motion that should be integrated [24]. The differential-algebraic equations governing a dynamical system contain two parts: the ODEs, which are the structural subsystem description, and the algebraic equations, which are the system's constraint equations. The constraints of a mechanical system are usually in a non-integrable velocity form, but many methods consider these equations in acceleration form [25]. That is to say, almost all traditional methods use constraint equations as a set of second-order equations along with second-order ODEs that are based on generalized coordinates and generalized velocities as motion variables. This is the first reason for the violation of system constraints and error accumulation. The second significant factor that contributes to the constraint (and other system's constants) drift is the number and complexity of the ODEs in DAEs (differential-algebraic equations). This is also the main determinant of CPU time.

In order to deal with these issues, many attempts were made [26]. One breakthrough was achieved by Gibbs-Appell by introducing the quasi-velocity notion [27]. The quasi-velocities were introduced as a substitution of generalized velocities to reduce the complexity of the equations, as well as their order [9]. The concept of quasi-velocities was further studied and used to eliminate Lagrange multipliers from motion equations in [28, 35]. The Gibbs-Appell method is still in the spotlight and has many subsequent variations. For instance, in [29], based on Gibss-Appell equations, an efficient formulation of motion equations is presented, in which there is no need to calculate Gibbs function. In [30] a concise and straightforward derivation of Gibbs-Appel equations of motion is presented. This type of system representation had a variety of benefits. Firstly, due to the flexibility of choosing the quasi-velocities, one can derive equations more easily and in a simpler form by correct selection of the quasi-velocities. Moreover, these dynamics equations use constraints equations in a velocity form, which leads to the reduction of the constraints violation error [28]. Besides, these constraint equations are embedded with the dynamical equation, diminishing the computational error. The mentioned advantages caused further interest in advancing these types of approaches [31-35] as well as increasing their application [36-40]. For instance, in [34], a canonical form of Kane method is developed for discrete dynamical systems. In [39], a snake robot is analyzed both in simulation and experiment using Gibbs-Appell approach. In [40], an efficient dynamic model for five-axis machine tools is provided in which there is no need to compute Lagrange multipliers. The method is proved to be 90% faster than existing methods. Despite order reduction and simplification of the equations, Kane and Gibbs-Appell-based approaches fail to give the minimum number of equations needed for analyzing a system including ignorable coordinates.

As mentioned earlier, Routh successfully tried to reduce the number of equations of an unconstrained system with ignorable coordinates. Assuming that the system has $m$ number of generalized coordinates, $p$ degrees of freedom, $r$ constraints and $s$ ignorable coordinates, Routh derived $m - s$ number of first-order equations. Using the Lagrange approach, in [41], a modified form of Routh equations, which is based on generalized coordinates and generalized velocities, was derived for both constrained and unconstrained dynamical systems. In that paper, a definition of the ignorable coordinates was introduced, and accordingly, a set of $p - s$ equations was derived in matrix form. In that research, the selection of the non-ignorable quasi-velocities was considered as a combination of only non-ignorable generalized velocities.

In this paper, an efficient form of Volterra's equations of motion for systems including ignorable coordinates is derived. This method is proved to be applicable for both constrained and unconstrained systems. In this approach, similar to the standard Volterra method, the quasi-velocity notation is used. Before derivation of the equations, based on the definition of the ignorable coordinates provided in [41], the notion of "dynamical constraints" is introduced. These constraints, along with ordinary (kinematical) constraints, are then embedded with Volterra's equations. The eventual equations are originally in the first-order form, and thus less effort is required to derive the efficient equations of motion. More importantly, in contrast to [41], the choice of non-ignorable quasi velocities is not restricted to non-ignorable generalized velocities, which means one can choose a combination of both ignorable and non-ignorable generalized velocities as non-ignorable quasi-velocities. This flexibility, if used properly, will lead to simpler equations and less computational effort. In the end, this approach derives a set of $p - s$ equations that are more satisfactory in terms of runtime and computational error than conventional methods.

The remainder of this paper is organized as follows: the proposed approach, which is the main contribution of this work, is provided in section 2. Three subsections are given in this section, that is, subsection 2.1., which provides the prerequisites and standard Volterra approach, subsection 2.2., which provides the reader with the concept and formulation of dynamical constraints, and subsection 2.3., which deals with the modified form of Volterra's equation. Moreover, a matrix notation of the proposed method is obtained, and a comparison (in terms of the number of equations) among conventional methods is drawn in this section. Section 3 presents the results of three illustrative simulation examples and a brief discussion about them. In section 4, conclusions are presented. The nomenclature list is provided in section 5.

## 2. The proposed method

### *2.1. Preliminaries and the standard Volterra approach*

Consider a dynamical system with $\boldsymbol{q} = [q_1 \quad \cdots \quad q_m]^T$ as its generalized coordinates vector. Suppose that the system has $r$ number of nonholonomic constraints as:

$$\boldsymbol{a}(t;\boldsymbol{q})\dot{\boldsymbol{q}} + \boldsymbol{b}(t;\boldsymbol{q}) = \boldsymbol{0} \tag{1}$$

in which $\boldsymbol{a} \in \mathbb{R}^r \times \mathbb{R}^m$ is constraints' Jacobian matrix, $\boldsymbol{b} \in \mathbb{R}^r$ is constraints' bias vector and $\dot{\boldsymbol{q}} \in \mathbb{R}^m$ is the generalized velocities vector. For a system with holonomic constraints (i.e., constraints that are defined on configuration) it is straightforward to achieve the form (1) by simply differentiating the holonomic form.

A special form of the fundamental equations of motion for a general multibody system is called Volterra's equation of motion, which for a system with $m$ GCs (generalized coordinates) and $r$ constraints (and therefore, $p = m - r$ Degrees of freedom (DOFs)) is formulated as [10]:

$$\frac{d}{dt}\left(\frac{\partial T}{\partial u_j}\right) - \sum_{i=1}^{N}\left[{\boldsymbol{P}_i}^T \frac{d}{dt}\left(\frac{\partial \boldsymbol{V}_{Gi}}{\partial u_j}\right) + {\boldsymbol{H}_{Gi}}^T \frac{d}{dt}\left(\frac{\partial \boldsymbol{\omega}_i}{\partial u_j}\right)\right] = U_j, \quad j = 1, \dots, p \tag{2}$$

in which $N$ is the number of the bodies of the system. Furthermore, $T$ is the kinetic energy of the system, $\boldsymbol{P}_i \in \mathbb{R}^3$ and $\boldsymbol{H}_{Gi} \in \mathbb{R}^3$ are the linear momentum and the angular momentum vectors of the $i$-th body of the system, $\boldsymbol{V}_{Gi} \in \mathbb{R}^3$ and $\boldsymbol{\omega}_i \in \mathbb{R}^3$ are the linear and the angular velocities of the center of the mass of the $i$-th body of the system, $\boldsymbol{u} \in \mathbb{R}^p$ is the quasi-velocity vector and $\boldsymbol{U} \in \mathbb{R}^p$ is the generalized forces vector expressed in quasi-velocity space.

The procedure of obtaining the equations of motion from (2) is as follows: first of all, one should choose $p$ number of quasi-velocities:

$$\boldsymbol{u} = \boldsymbol{Y}(t;\boldsymbol{q})\dot{\boldsymbol{q}} + \boldsymbol{Z}(t;\boldsymbol{q}) \tag{3}$$

in which $\boldsymbol{u} \in \mathbb{R}^p$, $\dot{\boldsymbol{q}} \in \mathbb{R}^m$, $\boldsymbol{Z} \in \mathbb{R}^p$ and $\boldsymbol{Y} \in \mathbb{R}^p \times \mathbb{R}^m$ are quasi-velocities, generalized velocities, bias vectors, and Jacobian matrix of quasi-velocities, respectively. Augmenting Eq. (3) with Eq. (1), one has:

$$\dot{\boldsymbol{q}} = \boldsymbol{W}(t;\boldsymbol{q})\boldsymbol{u} + \boldsymbol{X}(t;\boldsymbol{q}) \tag{4}$$

in which

$$\begin{cases} \boldsymbol{W} = \begin{bmatrix} \boldsymbol{Y} \\ \boldsymbol{a} \end{bmatrix}^{-1} \begin{bmatrix} \boldsymbol{I} \\ \boldsymbol{0} \end{bmatrix} \\ \boldsymbol{X} = -\begin{bmatrix} \boldsymbol{Y} \\ \boldsymbol{a} \end{bmatrix}^{-1} \begin{bmatrix} \boldsymbol{Z} \\ \boldsymbol{b} \end{bmatrix} \end{cases} \tag{5}$$

As the second step, the angular and linear velocities of the bodies of the system should be derived. It should be noted that originally these velocities are functions of time, generalized coordinates, and generalized velocities $(t, \boldsymbol{q}, \dot{\boldsymbol{q}})$, however, by using Eq. (4), one can transform these velocities into the form:

$$\begin{cases} \boldsymbol{V}_{Gi} = \boldsymbol{V}_{Gi}(t; \boldsymbol{q}; \boldsymbol{u}) \\ \boldsymbol{\omega}_i = \boldsymbol{\omega}_i(t; \boldsymbol{q}; \boldsymbol{u}) \end{cases} \tag{6}$$

Now, it's easy to calculate the kinetic energy:

$$T = \sum_{i=1}^{N} \frac{1}{2} (m_i \boldsymbol{V}_{Gi} . \boldsymbol{V}_{Gi} + \boldsymbol{\omega}_i . \boldsymbol{I}_{Gi} \boldsymbol{\omega}_i) \tag{7}$$

in which $m_i$ and $\boldsymbol{I}_{Gi} \in \mathbb{R}^3 \times \mathbb{R}^3$ are the mass and centroidal moment of inertia of the $i$-th body of the system, respectively. The next step is to calculate linear and angular momentums:

$$\begin{cases} \boldsymbol{P}_i = m_i \boldsymbol{V}_{Gi} \\ \boldsymbol{H}_{Gi} = \boldsymbol{I}_{Gi} \boldsymbol{\omega}_i \end{cases} \tag{8}$$

Then, the vector of all non-conservative and conservative generalized forces in the quasi-velocity space is to be determined as:

$$U_j = \sum_{i=1}^{\mu} \boldsymbol{F}_i . \left( \frac{\partial \boldsymbol{V}_{Gi}}{\partial u_j} \right) + \sum_{i=1}^{\lambda} \boldsymbol{M}_i . \left( \frac{\partial \boldsymbol{\omega}_i}{\partial u_j} \right) \tag{9}$$

in which $\boldsymbol{F}_i$ is the $i$-th external force, $\boldsymbol{M}_i$ is the $i$-th torque, $\mu$ is the number of forces applied on the system, and $\lambda$ is the number of torques acting on the system.

Lastly, using Eq. (2) one can simply derive the equations of the motion.

### *2.2. Formulating dynamical constraints*

Constraints of the form (1) are imposed on a system due to the presence of some limiting elements or some inherent restrictions that come from the nature of the system. For instance, when a wheel is mounted on a body, assuming that there is no slip between the wheel and the underlying ground, a velocity constraint is imposed on the system (which is the case in the first simulation example in section 3.1). Another example of this constraint is when a particular body of a system is required to move along a specific path, which is common especially in linkage mechanisms. These constraints arise from the physical limitations, geometry, or configuration of the system and hence, can be deemed as kinematical constraints. For systems including ignorable coordinates, another set of constraints can be defined. Equations of these constraints are derived from the definition of the ignorable coordinates. The derivation procedure is to be explained hereunder.

Suppose that a multibody system has $s$ number of ignorable coordinates (which will be defined later on), and therefore, one can split the vector of generalized coordinates into two parts:

$$\boldsymbol{q} = [\boldsymbol{q}_{NI}{}^T \quad , \quad \boldsymbol{q}_I{}^T]^T \tag{10}$$

in which $\boldsymbol{q}_{NI} \in \mathbb{R}^{(m-s)}$ and $\boldsymbol{q}_I \in \mathbb{R}^s$ are non-ignorable and ignorable GCs respectively. The same division can be made for the constraints' Jacobian and the vector of generalized forces:

$$\boldsymbol{a} = [\boldsymbol{a}_1 \quad , \quad \boldsymbol{a}_2] \tag{11}$$

$$\boldsymbol{Q} = [\boldsymbol{Q}_{NI}{}^T \quad , \quad \boldsymbol{Q}_I{}^T]^T \tag{12}$$

where $\boldsymbol{a}_1 \in \mathbb{R}^r \times \mathbb{R}^{(m-s)}$, $\boldsymbol{a}_2 \in \mathbb{R}^r \times \mathbb{R}^s$, $\boldsymbol{Q}_{NI} \in \mathbb{R}^{(m-s)}$ and $\boldsymbol{Q}_I \in \mathbb{R}^s$. Note that $\boldsymbol{Q} \in \mathbb{R}^m$ is the generalized forces vector expressed in the generalized velocity space, not in quasi-velocity space.

Based on [41], ignorable coordinates can be defined as follows:

*Definition1*. If the following conditions hold:

- The last $s$ rows of $\boldsymbol{q}$ in Eq. (10) do not appear in Lagrangian function, which by definition is $\mathcal{L} = T - V$, where $T$ and $V$ are kinetic and potential energies of the system, respectively.
- The last s columns of $\mathbf{a}$ in Eq. (11) are identically zero (i.e., $\mathbf{a}_2 = \mathbf{0}$).
- The last s components of $\mathbf{Q}$ in Eq. (12) are identically zero (i.e., $\mathbf{Q}_\mathrm{I} = \mathbf{0}$).

Then, the $\boldsymbol{q}_I$ is said to be the ignorable coordinates vector. Practical examples of ignorable coordinates are provided in section 3.

Eq. (7) can be rewritten in generalized velocity space as:

$$T(t;\boldsymbol{q};\dot{\boldsymbol{q}}) = \frac{1}{2}\dot{\boldsymbol{q}}^T\boldsymbol{M}(t;\boldsymbol{q}_{NI})\dot{\boldsymbol{q}} + \dot{\boldsymbol{q}}^T\boldsymbol{N}(t;\boldsymbol{q}_{NI}) + T_0(t;\boldsymbol{q}_{NI}) \tag{13}$$

in which $\boldsymbol{M} \in \mathbb{R}^m \times \mathbb{R}^m$, $\boldsymbol{N} \in \mathbb{R}^m$ and $T_0 \in \mathbb{R}$. It is important to note that $\boldsymbol{M}$ is the system's mass matrix, which is a symmetric positive definite matrix. Matrices $\boldsymbol{M}$ and $\boldsymbol{N}$ can be partitioned as:

$$\boldsymbol{M} = \begin{bmatrix} \boldsymbol{M}_{11} & \boldsymbol{M}_{12} \\ \boldsymbol{M}_{21} & \boldsymbol{M}_{22} \end{bmatrix} \tag{14}$$

$$\boldsymbol{N} = [\boldsymbol{N}_1{}^T \quad \boldsymbol{N}_2{}^T] \tag{15}$$

where $\boldsymbol{M}_{11}(t;\boldsymbol{q}_{NI}) \in \mathbb{R}^{(m-s)} \times \mathbb{R}^{(m-s)}$, $\boldsymbol{M}_{12}(t;\boldsymbol{q}_{NI}) \in \mathbb{R}^{(m-s)} \times \mathbb{R}^{s}$, $\boldsymbol{M}_{21}(t;\boldsymbol{q}_{NI}) \in \mathbb{R}^{s} \times \mathbb{R}^{(m-s)}$, $\boldsymbol{M}_{22}(t;\boldsymbol{q}_{NI}) \in \mathbb{R}^{s} \times \mathbb{R}^{s}$, $\boldsymbol{N}_1(t;\boldsymbol{q}_{NI}) \in \mathbb{R}^{(m-s)}$ and $\boldsymbol{N}_2(t;\boldsymbol{q}_{NI}) \in \mathbb{R}^{s}$. It should be noted that $\boldsymbol{M}_{21} = \boldsymbol{M}_{12}{}^T$ due to symmetricity of the $\boldsymbol{M}$.

Using Eq. (14) and Eq. (15), Eq. (13) can be rewritten as:

$$T(t,\boldsymbol{q},\dot{\boldsymbol{q}}) = \frac{1}{2}\dot{\boldsymbol{q}}_{NI}{}^T\boldsymbol{M}_{11}\dot{\boldsymbol{q}}_{NI} + \dot{\boldsymbol{q}}_{NI}{}^T\boldsymbol{M}_{12}\dot{\boldsymbol{q}}_I + \frac{1}{2}\dot{\boldsymbol{q}}_I{}^T\boldsymbol{M}_{22}\dot{\boldsymbol{q}}_I + \dot{\boldsymbol{q}}_{NI}{}^T\boldsymbol{N}_1 + \dot{\boldsymbol{q}}_I{}^T\boldsymbol{N}_2 + T_0 \quad (16)$$

in which $\dot{\boldsymbol{q}}_{NI} \in \mathbb{R}^{(m-s)}$ and $\dot{\boldsymbol{q}}_I \in \mathbb{R}^{s}$ are non-ignorable and ignorable generalized velocities vectors, respectively. Lagrange equations for a constrained dynamical system can be written as:

$$\frac{d}{dt}\left(\frac{\partial \mathcal{L}}{\partial \dot{\boldsymbol{q}}_j}\right) - \left(\frac{\partial \mathcal{L}}{\partial \boldsymbol{q}_j}\right) = \boldsymbol{Q}_j + \boldsymbol{a}^{\boldsymbol{T}}\boldsymbol{\lambda}, \quad j = 1,\dots,m \quad (17)$$

in which $\boldsymbol{\lambda} \in \mathbb{R}^{r}$ is the Lagrange multipliers vector. The Lagrange equations can be written for ignorable and non-ignorable coordinates, separately

$$\frac{d}{dt}\left(\frac{\partial \mathcal{L}}{\partial \dot{\boldsymbol{q}}_{NI}}\right) - \left(\frac{\partial \mathcal{L}}{\partial \boldsymbol{q}_{NI}}\right) = \boldsymbol{Q}_{NI} + \boldsymbol{a}_1{}^{\boldsymbol{T}}\boldsymbol{\lambda}_{NI} \quad (18)$$

$$\frac{d}{dt}\left(\frac{\partial \mathcal{L}}{\partial \dot{\boldsymbol{q}}_{I}}\right) - \left(\frac{\partial \mathcal{L}}{\partial \boldsymbol{q}_{I}}\right) = \boldsymbol{Q}_{I} + \boldsymbol{a}_2{}^{\boldsymbol{T}}\boldsymbol{\lambda}_{I} \quad (19)$$

By employing the definition of the ignorable coordinates provided in *Definition1*, Eq. (19) can be simplified to:

$$\frac{d}{dt}\left(\frac{\partial \mathcal{L}}{\partial \dot{\boldsymbol{q}}_{I}}\right) = 0 \quad (20)$$

or:

$$\frac{\partial \mathcal{L}}{\partial \dot{\boldsymbol{q}}_{I}} = \boldsymbol{G}_I = constant \quad (21)$$

in which $\boldsymbol{G}_I$ is a constant vector and is referred to as generalized momentum vector corresponding to the ignorable coordinates. Since $\mathcal{L} = T - V$, and the potential energy is not function of the generalized velocities (i.e., $\dot{\boldsymbol{q}}_I$), Eq. (21) can be simplified to:

$$\frac{\partial T}{\partial \dot{\boldsymbol{q}}_{I}} = \boldsymbol{G}_I = constant \quad (22)$$

From Eq. (16) it's easy to calculate the left hand side of Eq. (22) as:

$$\frac{\partial T}{\partial \dot{\boldsymbol{q}}_I} = \boldsymbol{M}_{21}\dot{\boldsymbol{q}}_{NI} + \boldsymbol{M}_{22}\dot{\boldsymbol{q}}_I + \boldsymbol{N}_2 \tag{23}$$

or:

$$\frac{\partial T}{\partial \dot{\boldsymbol{q}}_I} = \boldsymbol{M}'\dot{\boldsymbol{q}} + \boldsymbol{N}_2 \tag{24}$$

in which $\boldsymbol{M}' = [\boldsymbol{M}_{21} \quad \boldsymbol{M}_{22}]$. Using Eq. (22) and Eq. (24) we have:

$$\boldsymbol{M}'\dot{\boldsymbol{q}} + \boldsymbol{N}' = \boldsymbol{0} \tag{25}$$

in which $\boldsymbol{N}' = \boldsymbol{N}_2 - \boldsymbol{G}_I$. Eq. (25) can be considered as the dynamical constraint of the system. This is because it has the regular form of constraints (same as Eq. (1)), but is associated with the dynamics of the system. In Eq. (25), $\boldsymbol{M}' \in \mathbb{R}^s \times \mathbb{R}^m$ is constraint's Jacobian matrix and $\boldsymbol{N}' \in \mathbb{R}^s$ is constraint's bias vector.

The constant $\boldsymbol{G}_I$ can be easily computed using initial condition:

$$\boldsymbol{G}_I = \boldsymbol{M}_{21}\dot{\boldsymbol{q}}_{NI}(t=0) + \boldsymbol{M}_{22}\dot{\boldsymbol{q}}_I(t=0) + \boldsymbol{N}_2(t=0) \tag{26}$$

### *2.3. The modified form of Volterra's equations*

For a system with $m$ generalized coordinates, $s$ ignorable coordinates, $p$ degrees of freedom, and $r$ number of kinematic constraints, the minimum number of motion equations needed is $p - s$. However, as stated earlier, the standard Volterra method gives us $p$ number of equations. We now prove that the Volterra equation provided in Eq. (2), can be written for $p - s$ number of quasi-velocities instead of for $p$ number of them. For a system consisting of $N$ particles, assuming that the D'Alembert's principle applies, one can write:

$$\sum_{i=1}^{N} (\boldsymbol{F}_i - m_i \boldsymbol{a}_i).\delta \boldsymbol{r}_i = 0 \tag{27}$$

in which $\boldsymbol{F}_i$ is the resultant of forces applied on the $i$-th particle, $m_i$ is the mass of the $i$-th particle, $\boldsymbol{a}_i$ is its acceleration and $\delta \boldsymbol{r}_i$ is its virtual displacement vector. We choose $p - s$ number of quasi-velocities (equal to the minimum number of equations needed for motion analysis of the system). These quasi-velocities are regarded as non-ignorable quasi-velocities. The velocity of the $i$-th particle can be written as:

$$\boldsymbol{V}_i = \sum_{j=1}^{p-s} \boldsymbol{V}_{ij} u_{NI_j} + \boldsymbol{V}_{it} \tag{28}$$

in which $\boldsymbol{V}_{ij}$ is the partial velocity of the $i$-th particle and is defined as:

$$V_{ij} = \frac{\partial \boldsymbol{V}_i}{\partial u_{NI_j}}, \qquad i = 1, \dots, N, \qquad j = 1, \dots, p - s \tag{29}$$

$\boldsymbol{V}_{it}$ can also be defined accordingly. The differential form of Eq. (28) is:

$$d\boldsymbol{r}_i = \sum_{j=1}^{p-s} \boldsymbol{V}_{ij} d\gamma_{NI_j} + \boldsymbol{V}_{it} dt \tag{30}$$

in which $\gamma_j$ is the $j$-th quasi-coordinate. The variational form of Eq. (30) can be written as:

$$\delta\boldsymbol{r}_i = \sum_{j=1}^{p-s} \boldsymbol{V}_{ij} \delta\gamma_{NI_j} = \sum_{j=1}^{p-s} \frac{\partial \boldsymbol{V}_i}{\partial u_j} \delta\gamma_{NI_j} \tag{31}$$

Substituting Eq. (31) in Eq. (27) gives us:

$$\sum_{j=1}^{p-s} \left( \sum_{i=1}^{N} (\boldsymbol{F}_i - m_i \boldsymbol{a}_i) . \frac{\partial \boldsymbol{V}_i}{\partial u_{NI_j}} \right) \delta\gamma_{NI_j} = 0 \tag{32}$$

From Eq. (9), it is known that the generalized forces vector in the quasi-velocity space can be calculated through:

$$U_{NI_j} = \sum_{i=1}^{N} \boldsymbol{F}_i . \frac{\partial \boldsymbol{V}_i}{\partial u_{NI_j}} \tag{33}$$

Using Eq. (33), Eq. (32) can be rewritten as:

$$\sum_{j=1}^{p-s} \left( \sum_{i=1}^{N} -m_i \boldsymbol{a}_i . \frac{\partial \boldsymbol{V}_i}{\partial u_{NI_j}} + U_{NI_j} \right) \delta\gamma_{NI_j} = 0 \tag{34}$$

Since quasi-coordinates are independent, one can conclude:

$$\sum_{i=1}^{N} -m_i \boldsymbol{a}_i . \frac{\partial V_i}{\partial u_{NI_j}} + U_{NI_j} = 0 \tag{35}$$

Similar to Eq. (7), the kinetic energy of the system of particles can be rewritten as:

$$T = \sum_{i=1}^{N} \frac{1}{2} (m_i \boldsymbol{V}_i . \boldsymbol{V}_i) \tag{36}$$

Using Eq. (36), we can calculate:

$$\frac{\partial T}{\partial u_{NI_j}} = \sum_{i=1}^{N}\left(m_i \boldsymbol{V}_i . \frac{\partial \boldsymbol{V}_i}{\partial u_{NI_j}}\right) \tag{37}$$

Differentiating Eq. (37) with respect to time gives us:

$$\frac{d}{dt}\left(\frac{\partial T}{\partial u_{NI_j}}\right) = \sum_{i=1}^{N}\left(m_i \boldsymbol{a}_i . \frac{\partial \boldsymbol{V}_i}{\partial u_{NI_j}} + m_i \boldsymbol{V}_i . \frac{d}{dt}\left(\frac{\partial \boldsymbol{V}_i}{\partial u_{NI_j}}\right)\right) \tag{38}$$

or:

$$\sum_{i=1}^{N} m_i \boldsymbol{a}_i . \frac{\partial \boldsymbol{V}_i}{\partial u_{NI_j}} = \frac{d}{dt}\left(\frac{\partial T}{\partial \boldsymbol{u}_{NI_j}}\right) - \sum_{i=1}^{N} m_i \boldsymbol{V}_i . \frac{d}{dt}\left(\frac{\partial \boldsymbol{V}_i}{\partial u_{NI_j}}\right) \tag{39}$$

Using Eq. (39) and Eq. (35), we have:

$$\frac{d}{dt}\left(\frac{\partial T}{\partial u_{NI_j}}\right) - \sum_{i=1}^{N} \boldsymbol{P}_i . \frac{d}{dt}\left(\frac{\partial \boldsymbol{V}_i}{\partial u_{NI_j}}\right) = U_{NI_j}, \quad j = 1, \dots, p-s \tag{40}$$

Eq. (40) is the reduced form of Volterra's equations of motion for a system including $N$ particles. It is straightforward to obtain the same equation for a multibody system including rigid bodies. Via the same procedure, we have:

$$\frac{d}{dt}\left(\frac{\partial T}{\partial u_{NI_j}}\right) - \sum_{i=1}^{N}\left[\boldsymbol{P}_i . \frac{d}{dt}\left(\frac{\partial \boldsymbol{V}_{Gi}}{\partial u_{NI_j}}\right) + \boldsymbol{H}_{Gi} . \frac{d}{dt}\left(\frac{\partial \boldsymbol{\omega}_i}{\partial u_{NI_j}}\right)\right] = U_{NI_j}, \quad j = 1, \dots, p-s \tag{41}$$

Eq. (41) is the reduced form of Volterra's equations of motion for a multibody system including ignorable coordinates, which leads to $p - s$ number of equations (i.e., the minimal number of equations).

In the remainder of this section, after obtaining the matrix notation of the proposed method, we will draw a comparison, in terms of the number of equations, among the proposed method and other conventional methods. Now the matrix form of the Eq. (41) is to be obtained. First of all, the Lagrangian function should be formed, and using the definition of the ignorable coordinates, the dynamical constraint's equation (Eq. (25)) should be derived. Then, $p - s$ number of quasi-velocities should be selected. As mentioned earlier, unlike [41], the choices in this approach are not limited to the non-ignorable generalized velocities:

$$\boldsymbol{u}_{NI} = \boldsymbol{Y}_{NI}(t; \boldsymbol{q})\dot{\boldsymbol{q}} + \boldsymbol{Z}_{NI}(t; \boldsymbol{q}) \tag{42}$$

in which $\boldsymbol{u}_{NI} \in \mathbb{R}^{p-s}$, $\dot{\boldsymbol{q}} \in \mathbb{R}^m$, $\boldsymbol{Z}_{NI} \in \mathbb{R}^{p-s}$ and $\boldsymbol{Y}_{NI} \in \mathbb{R}^{p-s} \times \mathbb{R}^m$ are non-ignorable quasi-velocities, generalized velocities, bias vectors and Jacobian matrix, respectively. By augmenting Eq. (42) with kinematic constraints (Eq. (1)) and dynamical constraint (Eq. (25)), we have:

$$\dot{\boldsymbol{q}} = \boldsymbol{W}_{NI}(t;\boldsymbol{q})\boldsymbol{u}_{NI} + \boldsymbol{X}_{NI}(t;\boldsymbol{q}) \tag{43}$$

in which

$$\begin{cases} \boldsymbol{W}_{NI} = \begin{bmatrix} \boldsymbol{Y}_{NI} \\ \boldsymbol{M}' \\ \boldsymbol{a} \end{bmatrix}^{-1} \begin{bmatrix} \boldsymbol{I}_{(p-s)\times(p-s)} \\ \boldsymbol{0}_{s\times(p-s)} \\ \boldsymbol{0}_{r\times(p-s)} \end{bmatrix} \\ \boldsymbol{X}_{NI} = -\begin{bmatrix} \boldsymbol{Y}_{NI} \\ \boldsymbol{M}' \\ \boldsymbol{a} \end{bmatrix}^{-1} \begin{bmatrix} \boldsymbol{Z}_{NI} \\ \boldsymbol{N}' \\ \boldsymbol{b} \end{bmatrix} \end{cases} \tag{44}$$

Eq. (41) can be rewritten as:

$$\frac{d}{dt}\left(\frac{\partial T}{\partial \boldsymbol{u}_{NI}}\right) - \sum_{i=1}^{N}\left[\frac{d}{dt}\left(\frac{\partial \boldsymbol{V}_{Gi}}{\partial \boldsymbol{u}_{NI}}\right)^T \boldsymbol{P}_i + \frac{d}{dt}\left(\frac{\partial \boldsymbol{\omega}_i}{\partial \boldsymbol{u}_{NI}}\right)^T \boldsymbol{H}_{Gi}\right] = \boldsymbol{U}_{NI} \tag{45}$$

It is desired to obtain the matrix form of each term of Eq. (45). We start with the first expression in the left hand side of the equation. Using Eq. (13) and Eq. (43), we have:

$$T(t;\boldsymbol{q};\boldsymbol{u}_{NI}) = \frac{1}{2}\boldsymbol{u}_{NI}{}^T\boldsymbol{M}_{NI}(t;\boldsymbol{q})\boldsymbol{u}_{NI} + \boldsymbol{u}_{NI}{}^T\boldsymbol{N}_{NI}(t;\boldsymbol{q}) + T_{0,NI}(t;\boldsymbol{q}) \tag{46}$$

in which $\boldsymbol{M}_{NI} \in \mathbb{R}^{p-s} \times \mathbb{R}^{p-s}$, $\boldsymbol{N}_{NI} \in \mathbb{R}^{p-s}$, $T_{0,NI} \in \mathbb{R}$. Additionally:

$$\boldsymbol{M}_{NI} = \boldsymbol{W}_{NI}{}^T\boldsymbol{M}\boldsymbol{W}_{NI} \tag{47}$$

$$\boldsymbol{N}_{NI} = \boldsymbol{W}_{NI}{}^T(\boldsymbol{M}\boldsymbol{X}_{NI} + \boldsymbol{N}) \tag{48}$$

$$T_{0,NI} = \frac{1}{2}\boldsymbol{X}_{NI}{}^T\boldsymbol{M}\boldsymbol{X}_{NI} + \boldsymbol{X}_{NI}{}^T\boldsymbol{N} + T_0 \tag{49}$$

and other parameters are defined previously. Differentiating Eq. (46) with respect to $\boldsymbol{u}_{NI}$ gives:

$$\frac{\partial T}{\partial \boldsymbol{u}_{NI}} = \boldsymbol{M}_{NI}\boldsymbol{u}_{NI} + \boldsymbol{N}_{NI} \tag{50}$$

Likewise, by differentiating Eq. (50) with respect to time one has:

$$\frac{d}{dt}\left(\frac{\partial T}{\partial \boldsymbol{u}_{NI}}\right) = \boldsymbol{M}_{NI}\dot{\boldsymbol{u}}_{NI} + \boldsymbol{A}_{NI}(t;\boldsymbol{q};\boldsymbol{u}_{NI}) \tag{51}$$

in which $\boldsymbol{A}_{NI} \in \mathbb{R}^{p-s}$, and using Eq. (43) it is defined as:

$$\boldsymbol{A}_{NI} = \frac{\partial}{\partial t}\left(\frac{\partial T}{\partial \boldsymbol{u}_{NI}}\right) + \frac{\partial}{\partial \boldsymbol{q}}\left(\frac{\partial T}{\partial \boldsymbol{u}_{NI}}\right)(\boldsymbol{W}_{NI}\boldsymbol{u}_{NI} + \boldsymbol{X}_{NI}) \tag{52}$$

Now we move on to the terms that appeared in front of the sigma in Eq. (45) and try to extract a matrix form for them. Similar to what was mentioned earlier, linear and angular velocities of the bodies of the system including ignorable coordinates, are initially of the form:

$$\begin{cases} \boldsymbol{V}_{Gi} = \boldsymbol{V}_{Gi}(t; \boldsymbol{q}_{NI}; \dot{\boldsymbol{q}}) \\ \boldsymbol{\omega}_i = \boldsymbol{\omega}_i(t; \boldsymbol{q}_{NI}; \dot{\boldsymbol{q}}) \end{cases} \tag{53}$$

Eq. (53) can be written as:

$$\begin{cases} \boldsymbol{V}_{Gi} = \boldsymbol{B}_i(t; \boldsymbol{q}_{NI})\dot{\boldsymbol{q}} + \boldsymbol{C}_i(t; \boldsymbol{q}_{NI}) \\ \boldsymbol{\omega}_i = \boldsymbol{D}_i(t; \boldsymbol{q}_{NI})\dot{\boldsymbol{q}} + \boldsymbol{E}_i(t; \boldsymbol{q}_{NI}) \end{cases} \tag{54}$$

in which $\boldsymbol{B}_i \in \mathbb{R}^3 \times \mathbb{R}^m$ and $\boldsymbol{D}_i \in \mathbb{R}^3 \times \mathbb{R}^m$ can be termed as Jacobian matrices of linear and angular velocities of the *i*-th body of the system, respectively. Similarly, $\boldsymbol{C}_i \in \mathbb{R}^3$ and $\boldsymbol{E}_i \in \mathbb{R}^3$ can be called bias matrices of velocities of the bodies of the system.

By using Eq. (43) one can transform the Eq. (54) into the quasi-velocity space:

$$\begin{cases} \boldsymbol{V}_{Gi} = \boldsymbol{B}_i\boldsymbol{W}_{NI}\boldsymbol{u}_{NI} + \boldsymbol{B}_i\boldsymbol{X}_{NI} + \boldsymbol{C}_i \\ \boldsymbol{\omega}_i = \boldsymbol{D}_i\boldsymbol{W}_{NI}\boldsymbol{u}_{NI} + \boldsymbol{D}_i\boldsymbol{X}_{NI} + \boldsymbol{E}_i \end{cases} \tag{55}$$

Differentiating Eq. (55) with respect to $\boldsymbol{u}_{NI}$ leads to:

$$\begin{cases} \dfrac{\partial \boldsymbol{V}_{Gi}}{\partial \boldsymbol{u}_{NI}} = \boldsymbol{B}_i\boldsymbol{W}_{NI}(t; \boldsymbol{q}) \\ \dfrac{\partial \boldsymbol{\omega}_i}{\partial \boldsymbol{u}_{NI}} = \boldsymbol{D}_i\boldsymbol{W}_{NI}(t; \boldsymbol{q}) \end{cases} \tag{56}$$

Likewise, differentiating Eq. (56) with respect to time results in:

$$\begin{cases} \dfrac{d}{dt}\left(\dfrac{\partial \boldsymbol{V}_{Gi}}{\partial \boldsymbol{u}_{NI}}\right) = \boldsymbol{B}_{NI_i}(t; \boldsymbol{q}; \boldsymbol{u}_{NI}) \\ \dfrac{d}{dt}\left(\dfrac{\partial \boldsymbol{\omega}_i}{\partial \boldsymbol{u}_{NI}}\right) = \boldsymbol{C}_{NI_i}(t; \boldsymbol{q}; \boldsymbol{u}_{NI}) \end{cases} \tag{57}$$

in which $\boldsymbol{B}_{NI} \in \mathbb{R}^3 \times \mathbb{R}^{p-s}$ and $\boldsymbol{C}_{NI} \in \mathbb{R}^3 \times \mathbb{R}^{p-s}$. These matrices are calculated through:

$$\begin{cases} \boldsymbol{B}_{NI_i} = \dfrac{\partial}{\partial t}\left(\dfrac{\partial \boldsymbol{V}_{Gi}}{\partial \boldsymbol{u}_{NI}}\right) + \left[\dfrac{\partial}{\partial \boldsymbol{q}}\left(\dfrac{\partial \boldsymbol{V}_{Gi}}{\partial \boldsymbol{u}_{NI_1}}\right)(\boldsymbol{W}_{NI}\boldsymbol{u}_{NI} + \boldsymbol{X}_{NI}) \quad \dots \quad \dfrac{\partial}{\partial \boldsymbol{q}}\left(\dfrac{\partial \boldsymbol{V}_{Gi}}{\partial \boldsymbol{u}_{NI_{p-s}}}\right)(\boldsymbol{W}_{NI}\boldsymbol{u}_{NI} + \boldsymbol{X}_{NI})\right] \\ \boldsymbol{C}_{NI_i} = \dfrac{\partial}{\partial t}\left(\dfrac{\partial \boldsymbol{\omega}_i}{\partial \boldsymbol{u}_{NI}}\right) + \left[\dfrac{\partial}{\partial \boldsymbol{q}}\left(\dfrac{\partial \boldsymbol{\omega}_i}{\partial \boldsymbol{u}_{NI_1}}\right)(\boldsymbol{W}_{NI}\boldsymbol{u}_{NI} + \boldsymbol{X}_{NI}) \quad \dots \quad \dfrac{\partial}{\partial \boldsymbol{q}}\left(\dfrac{\partial \boldsymbol{\omega}_i}{\partial \boldsymbol{u}_{NI_{p-s}}}\right)(\boldsymbol{W}_{NI}\boldsymbol{u}_{NI} + \boldsymbol{X}_{NI})\right] \end{cases} \tag{58}$$

Using Eq. (8) and Eq. (55) it is straightforward to write:

$$\begin{cases} \boldsymbol{P}_i = m_i(\boldsymbol{B}_i\boldsymbol{W}_{NI}\boldsymbol{u}_{NI} + \boldsymbol{B}_i\boldsymbol{X}_{NI} + \boldsymbol{C}_i) \\ \boldsymbol{H}_{Gi} = \boldsymbol{I}_{Gi}(\boldsymbol{C}_i\boldsymbol{W}_{NI}\boldsymbol{u}_{NI} + \boldsymbol{C}_i\boldsymbol{X}_{NI} + \boldsymbol{E}_i) \end{cases} \tag{59}$$

By using Eq. (57) and Eq. (59), it is easy to obtain the second term in the left hand side of the Eq. (45) as:

$$\sum_{i=1}^{N}\left[\frac{d}{dt}\left(\frac{\partial \boldsymbol{V}_{Gi}}{\partial \boldsymbol{u}_{NI}}\right)^T \boldsymbol{P}_i + \frac{d}{dt}\left(\frac{\partial \boldsymbol{\omega}_i}{\partial \boldsymbol{u}_{NI}}\right)^T \boldsymbol{H}_{Gi}\right] = \sum_{i=1}^{N}\left(\boldsymbol{F}_{NI_i}(t,\boldsymbol{q},\boldsymbol{u}_{NI}) + \boldsymbol{J}_{NI_i}(t,\boldsymbol{q},\boldsymbol{u}_{NI})\right) \tag{60}$$

in which

$$\begin{cases} \boldsymbol{F}_{NI_i} = m_i\boldsymbol{B}_{NI_i}{}^{\boldsymbol{T}}(\boldsymbol{B}_i\boldsymbol{W}_{NI}\boldsymbol{u}_{NI} + \boldsymbol{B}_i\boldsymbol{X}_{NI} + \boldsymbol{C}_i) \\ \boldsymbol{J}_{NI_i} = \boldsymbol{C}_{NI_i}{}^{\boldsymbol{T}}\boldsymbol{I}_{Gi}(\boldsymbol{C}_i\boldsymbol{W}_{NI}\boldsymbol{u}_{NI} + \boldsymbol{C}_i\boldsymbol{X}_{NI} + \boldsymbol{E}_i) \end{cases} \tag{61}$$

Eq. (60) can be reworked as:

$$\sum_{i=1}^{N}\left[\frac{d}{dt}\left(\frac{\partial \boldsymbol{V}_{Gi}}{\partial \boldsymbol{u}_{NI}}\right)^T \boldsymbol{P}_i + \frac{d}{dt}\left(\frac{\partial \boldsymbol{\omega}_i}{\partial \boldsymbol{u}_{NI}}\right)^T \boldsymbol{H}_{Gi}\right] = \boldsymbol{K}_{NI} \tag{62}$$

in which $\boldsymbol{K}_{NI} \in \mathbb{R}^{p-s}$ , and is defined as:

$$\boldsymbol{K}_{NI}(t;\boldsymbol{q};\boldsymbol{u}_{NI}) = \sum_{i=1}^{N}\left(\boldsymbol{F}_{NI_i}(t,\boldsymbol{q},\boldsymbol{u}_{NI}) + \boldsymbol{J}_{NI_i}(t,\boldsymbol{q},\boldsymbol{u}_{NI})\right) \tag{63}$$

Using Eq. (62) and Eq. (51) in Eq. (45) gives us the matrix form of reduced Volterra's equations:

$$\boldsymbol{M}_{NI}(t;\boldsymbol{q})\dot{\boldsymbol{u}}_{NI} = \boldsymbol{L}_{NI}(t;\boldsymbol{q};\boldsymbol{u}_{NI}) \tag{64}$$

in which $\boldsymbol{L}_{NI} \in \mathbb{R}^{p-s}$, and is defined as:

$$\boldsymbol{L}_{NI} = \boldsymbol{U}_{NI} + \boldsymbol{F}_{NI} - \boldsymbol{A}_{NI} \tag{65}$$

In order to use the reduced form of Volterra's equations of motion, the state variables vector should be chosen as:

$$\boldsymbol{z}^T = [\boldsymbol{q}^T \quad , \quad \boldsymbol{u}_{NI}{}^T] \tag{66}$$

Then, using Eq. (43) and Eq. (64), the time derivative of the state variables vector is easily obtained:

$$\dot{\boldsymbol{z}} = \begin{Bmatrix} \boldsymbol{W}_{NI}(t;\boldsymbol{q})\boldsymbol{u}_{NI} + \boldsymbol{X}_{NI}(t;\boldsymbol{q}) \\ {\boldsymbol{M}_{NI}}^{-1}(t;\boldsymbol{q})\boldsymbol{L}_{NI}(t;\boldsymbol{q};\boldsymbol{u}_{NI}) \end{Bmatrix} \tag{67}$$

Eq. (67) is the state-space representation of the efficient form of Volterra's equations of motion for a multibody system including ignorable coordinates. As can be understood from Eq. (67), the matrix whose inverse must be calculated (i.e., the $\boldsymbol{M}_{NI}$) is a square matrix of the size $p - s$. Comparing Eq. (46) to Eq. (14), one can understand the resemblance between the mass matrix ($\boldsymbol{M}$) and this matrix ($\boldsymbol{M}_{NI}$). Therefore, we can call this matrix the reduced mass matrix of the system. Using the same procedure which was carried out to derive the state-space representation for the reduced Volterra method, it is straightforward to obtain the same representation for the standard Volterra method. After deriving such formulation, it could be easily seen that the matrix whose inverse should be computed is the system's mass matrix (i.e., $\boldsymbol{M}$). As a result, since in the standard Volterra method the inverse matrix of the size $m$ should be calculated, but in the reduced Volterra method the inverse of a $p - s$ by $p - s$ matrix should be computed, the computation time and error associated with the efficient Volterra method is expected to be less than the standard Volterra method.

Moreover, when compared to other conventional methods, the number of the equations in the proposed method (or equivalently, the size of the matrix whose inverse should be calculated) is minimum. More precisely, basic methods such as Lagrange and Maggi, which work in the generalized speed space, give $m$ number of equations. Although other common methods, like Kane (or Gibbs-Appell) which use the concept of quasi-velocity, omit $r$ unnecessary equations (associated with constraints) and thus reduce the number of equations to $p = m - r$, they do not produce the minimum number of equations, since they still include the $s$ redundant equations (associated with the ignorable coordinates). A comparison of some common analytical approaches for deriving motion equations is provided in Table 1.

It is also important to discuss the existence of $\boldsymbol{W}_{NI}$ which is defined in Eq. (44). It can be seen from the Eq. (44) that the formulation of the proposed method is based on the existence of $\boldsymbol{W}_{NI}$ and $\boldsymbol{X}_{NI}$. In other words, we shall prove that the following matrix is invertible:

$$\boldsymbol{A} = \begin{bmatrix} \boldsymbol{Y}_{NI} \\ \boldsymbol{M}' \\ \boldsymbol{a} \end{bmatrix}$$

Needless to say, one selects *independent* quasi-velocities, so $\boldsymbol{Y}_{NI}$ has $p - s$ independent rows. Additionally, since the system's mass matrix is positive definite, $\boldsymbol{M}'$ has $s$ independent rows. Moreover, we assume that the kinematic constraints are independent, so $\boldsymbol{a}$ has $r$ independent rows. As a result, the matrix $\boldsymbol{A}$ has $p - s + s + r = m$ independent rows and thus is invertible. Now, a question may arise: How can we guarantee that the rows of one matrix are independent of the rows of the other two matrices? For instance, are the rows of $\boldsymbol{Y}_{NI}$ independent of the rows of $\boldsymbol{M}'$ or $\boldsymbol{a}$? The answer lies in the fact that the sources of these three matrices are different. To be more specific, $\boldsymbol{Y}_{NI}$ is originated from our selection of quasi-velocities, $\boldsymbol{M'}$ is derived from the system's properties, and $\boldsymbol{a}$ is obtained from the kinematics. In real physical systems, these matrices are independent due to their different characteristics where one comes from dynamics and the other one comes from kinematics. If the quasi-velocities are not chosen properly, matrix $\boldsymbol{A}$ may be non-invertible. Such cases could be detected and avoided.

**Table 1.** Comparison of various analytical methods in terms of number of equations

| ***Method*** | Lagrange | Maggi | Udwiadia | Boltzmann-Hamel | Volterra | Kane (Gibbs-Appell) | Routh | Proposed (efficeint Volta) |
|---|---|---|---|---|---|---|---|---|
| ***Number of equations*** | $m$ | $m$ | $m$ | $p$ | $p$ | $p$ | $m-s$ | $p-s$ |

## 3. Simulation results and discussion

In this section three multibody systems which include ignorable coordinates are considered in order to analyze the performance of the proposed method and compare it to conventional methods of deriving the equations of motion.

### *3.1. A cart with a 2-DOF pendulum (constrained system)*

The system shown in Figure 1, consisting of a cart (with $G$ as its center of mass), a 2-DOF pendulum, and a motor mounted on point $B$, is to be considered as the first case study. The vector of generalized coordinates is chosen as:

$$\boldsymbol{q} = [\theta_1 \quad \theta_2 \quad x]^T \tag{68}$$

in which $x$ is measured from the left side of the cart at the rest point. As shown in Figure 1, there is a wheel mounted on point $C$ aligned with the bar $BC$. This wheel imposes a velocity constraint on the system, that is, the velocity of point $C$ relative to the cart should be along the direction of the bar $BC$. Mathematically, this constraint can be written as:

$$[l\,cos(\theta_1 - \theta_2) \quad l \quad 0]\begin{Bmatrix} \dot{\theta}_1 \\ \dot{\theta}_2 \\ \dot{x} \end{Bmatrix} = 0 \tag{69}$$

Therefore, the system has 2 DOFs. Moreover, one can calculate:

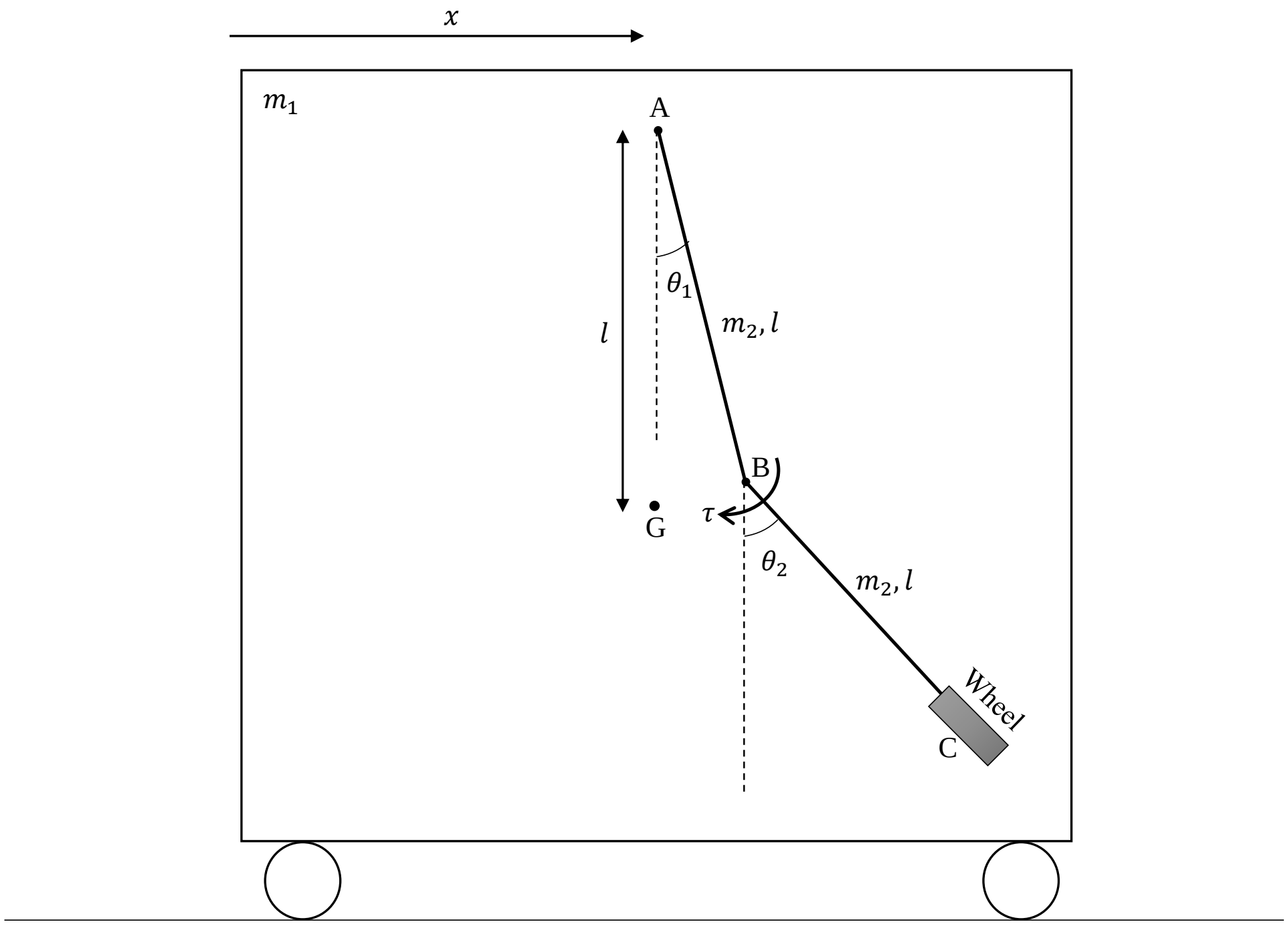

**Fig. 1** First simulation example: a cart with a 2-DOF pendulum

$$\boldsymbol{Q} = \begin{Bmatrix} \tau \\ -\tau \\ 0 \end{Bmatrix} \tag{70}$$

For methods which require $p$ number of quasi-velocities, the quasi-velocity vector is chosen as:

$$\boldsymbol{u} = \begin{Bmatrix} \dot{\theta}_1 - \dot{\theta}_2 \\ \dot{x} \end{Bmatrix} \tag{71}$$

After calculating the Lagrangian function, it can be easily seen that it is not a function of $x$. Moreover, comparing Eq. (69) and Eq. (70) to Eq. (11) and Eq. (12), it is easy to observe that $\boldsymbol{a}_2 = \boldsymbol{0}$ and $\boldsymbol{Q}_\mathrm{I} = \boldsymbol{0}$. As a result, according to *Definition1* it can be concluded that $x$ is an ignorable coordinate. The initial condition and properties of the system are provided in Table 2.

**Table 2.** System parameters and initial conditions (the first case study) [41]

| ***Parameter*** | ***Value*** | ***Unit*** |
|---|---|---|
| $m_1$ | 1 | $kg$ |
| $m_2$ | 0.5 | $kg$ |
| $l$ | 0.2 | $m$ |
| $\theta_1$ | $\pi/2$ | $rad$ |
| $\theta_2$ | $\pi/2$ | $rad$ |
| $x$ | 4 | $m$ |
| $\dot{\theta}_1$ | 1 | $rad/s$ |
| $\dot{\theta}_2$ | -1 | $rad/s$ |
| $\dot{x}$ | 3 | $m/s$ |

Moreover, the value of the torque ($\tau$) is set to zero. In this case, due to the lack of external active forces and torques, the amount of mechanical energy of the system must remain constant. Therefore, by computing the error of mechanical energy conservation for each method, in addition to linear momentum ($X$ direction) conservation error, an extra criterion for comparison of the methods is available.

The problem is solved using 4 different methods: Lagrange, Maggi, Gibbs-Appell (Kane), and the proposed method. Quasi-velocities in Maggi and Gibbs-Appell methods are chosen as Eq. (71), and in the proposed method, the non-ignorable quasi-velocity is chosen as $u_{NI} = \dot{\theta}_2 - \dot{\theta}_1$. The simulation is carried out for 50 seconds, the solving time step is fixed at 0.01 (s) and the ode45 function of MATLAB is used to solve the equations. The results of the simulation are provided in Table 3.

**Table 3.** Result of simulation using various methods (first case study)

| ***Method*** | Lagrange | Maggi | Gibbs-Appell (Kane) | Proposed |
|---|---|---|---|---|
| ***Number of State Variables*** | 6 | 6 | 5 | 4 |
| ***Number of Equations*** | 3 | 3 | 2 | 1 |
| ***CPU Time (s)*** | 1.93 | 1.72 | 1.87 | 1.69 |
| ***Norm of Energy Error*** | $6.05 \times 10^{-4}$ | $6.05 \times 10^{-4}$ | $9.77 \times 10^{-5}$ | $3.55 \times 10^{-6}$ |
| ***Norm of Constraint Error*** | $7.51 \times 10^{-5}$ | $7.52 \times 10^{-5}$ | $6.54 \times 10^{-15}$ | 0 |
| ***Norm of Linear Momentum Error (X)*** | $3.83 \times 10^{-6}$ | $3.93 \times 10^{-6}$ | $3.96 \times 10^{-5}$ | $2.41 \times 10^{-15}$ |

As it's understandable from Table 3, the proposed method's performance is considerably better than other approaches that were used to solve the problem. The plots of mechanical energy drift, constraint's conservation error, and linear momentum conservation error (in $X$ direction) are shown in Figures 2 to 4. Again, the preferable results of the proposed method can be concluded in comparison with alternative approaches. The plot of the generalized coordinates is illustrated in Figure 5. Additionally, shots of the motion of the system in several instants are displayed in Figure 6.

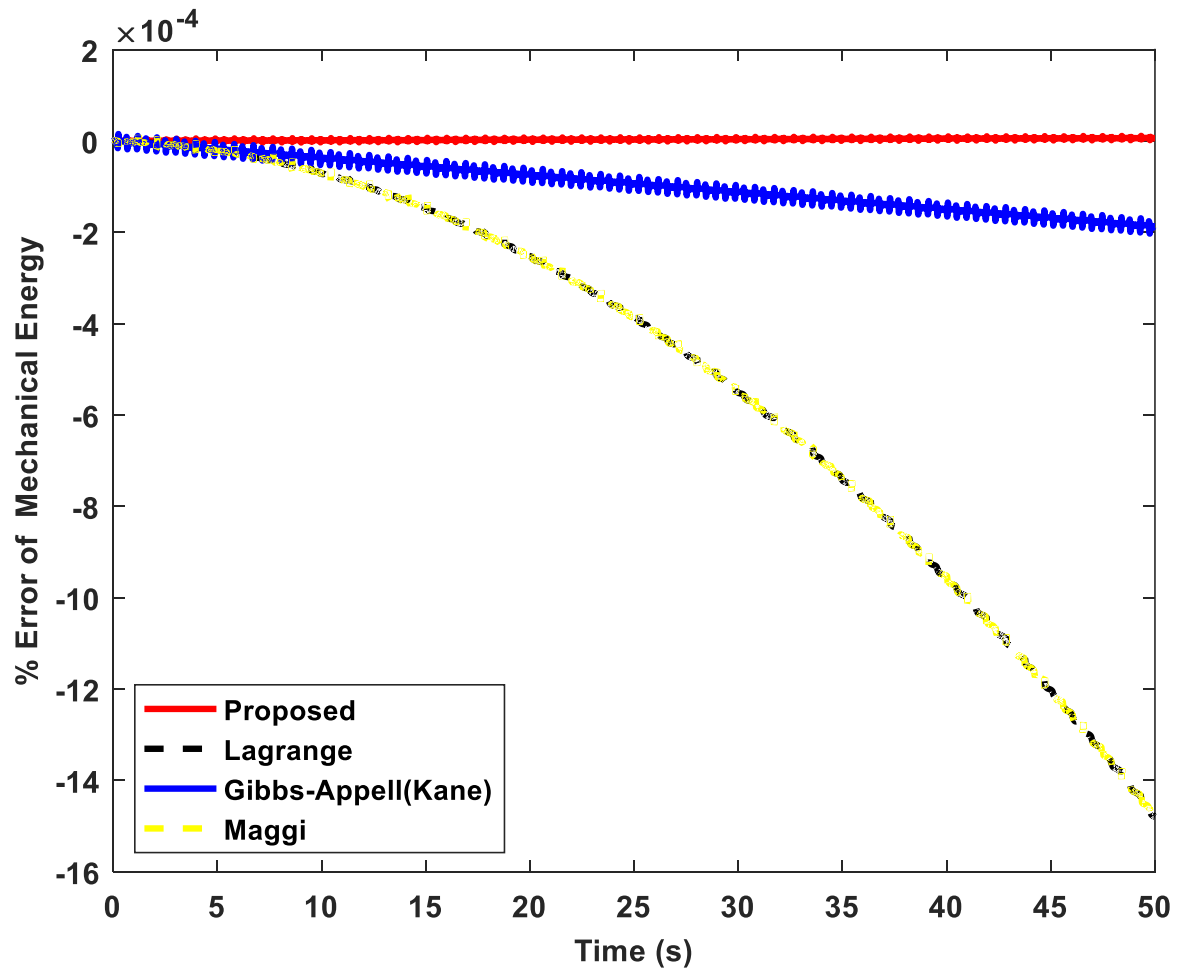

**Fig. 2** Percentage of mechanical energy drift in the first case study

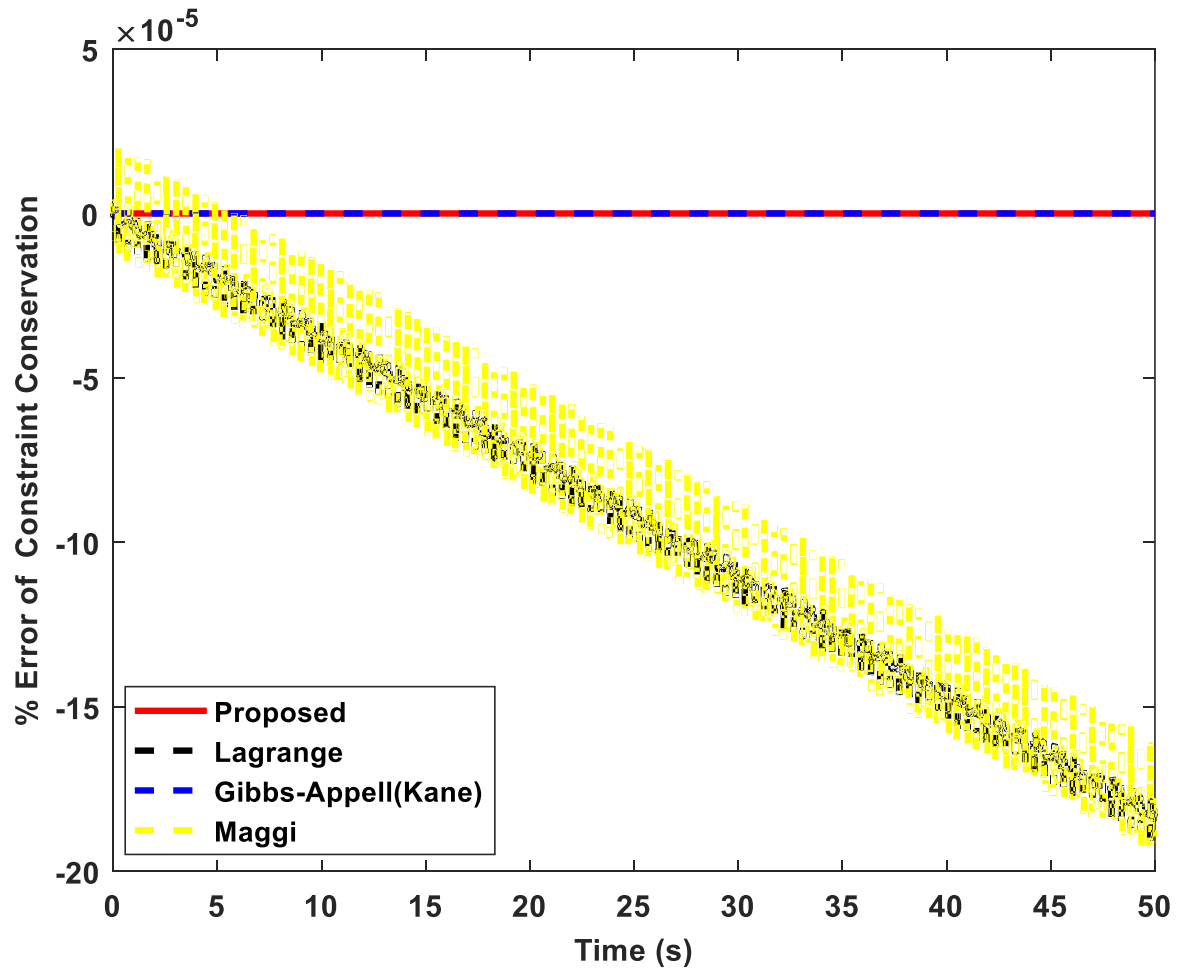

**Fig. 3** Percentage of constraint conservation error in the first case study

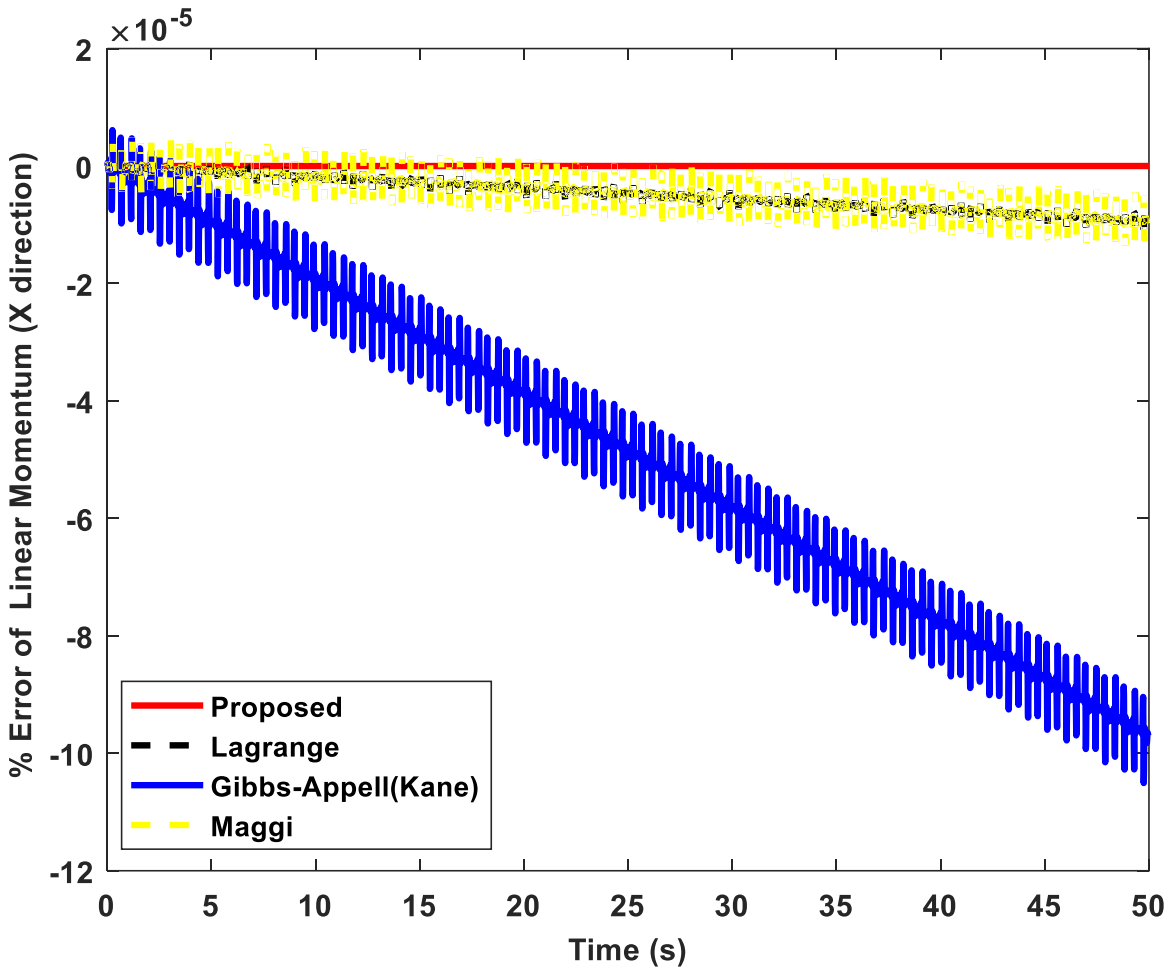

**Fig. 4** Percentage of linear momentum (x direction) conservation error in the first case study

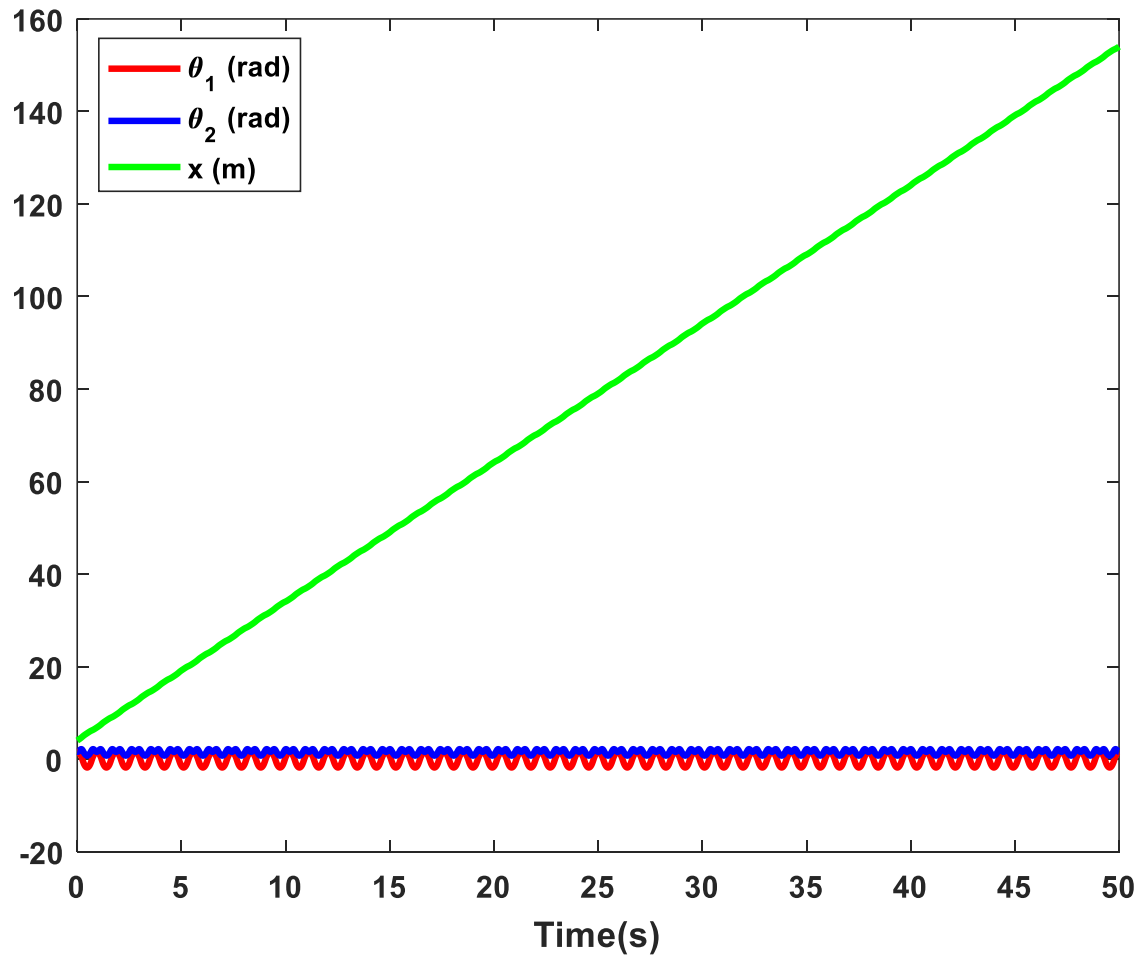

**Fig. 5** The plot of the generalized coordinates in the first case study

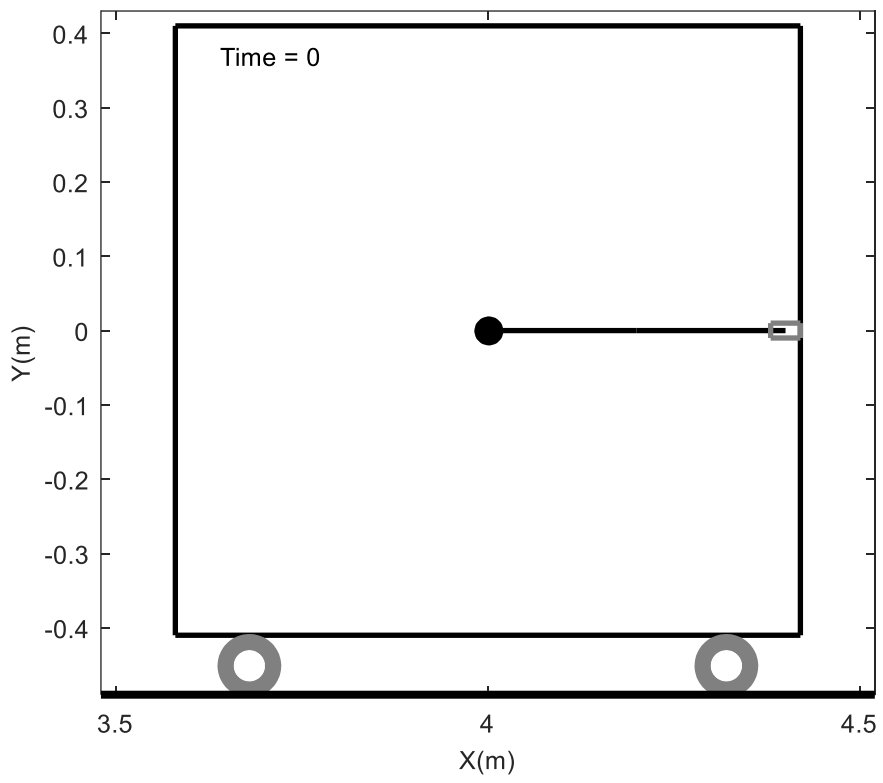

(a)

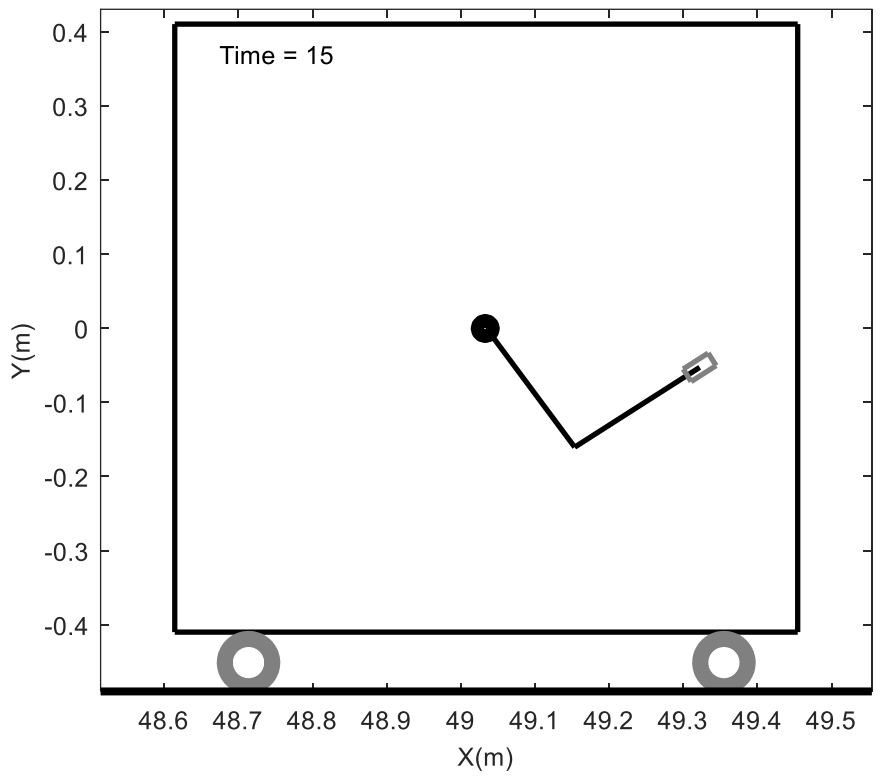

(b)

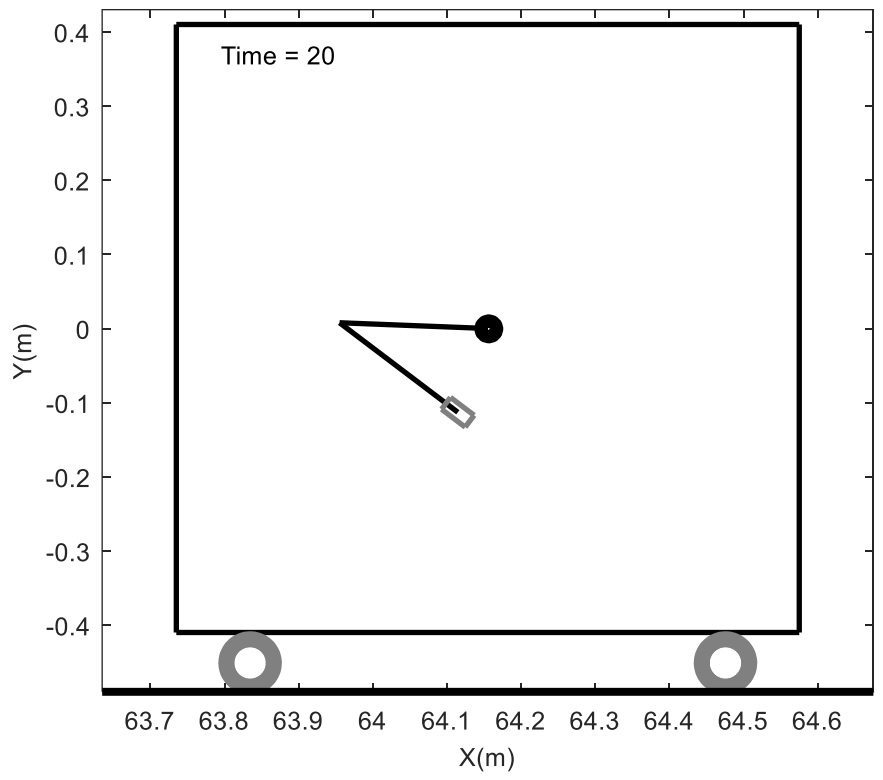

(c)

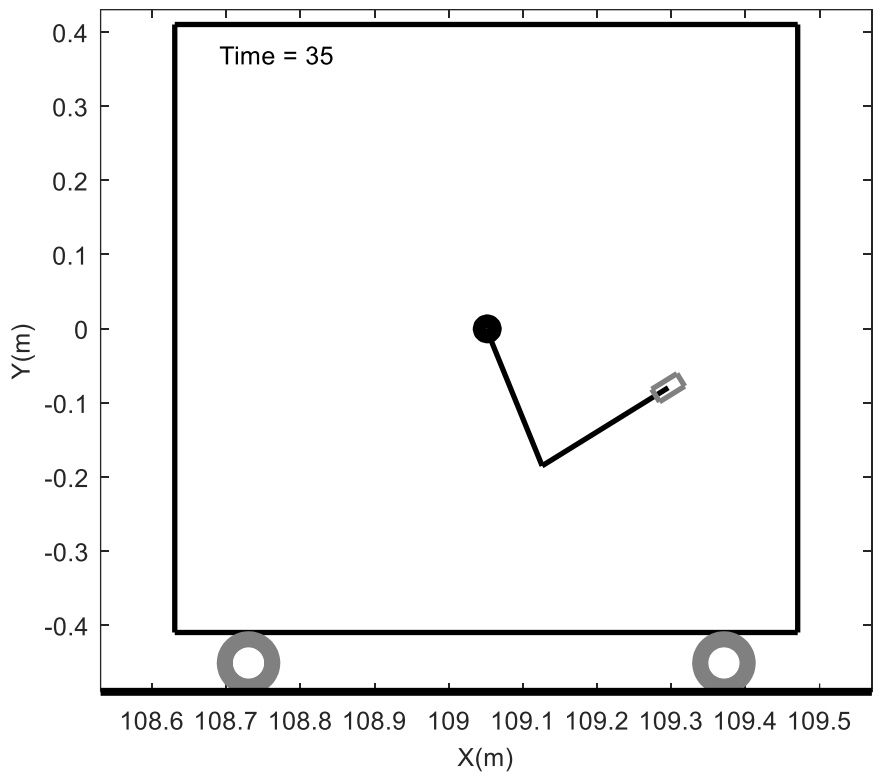

(d)

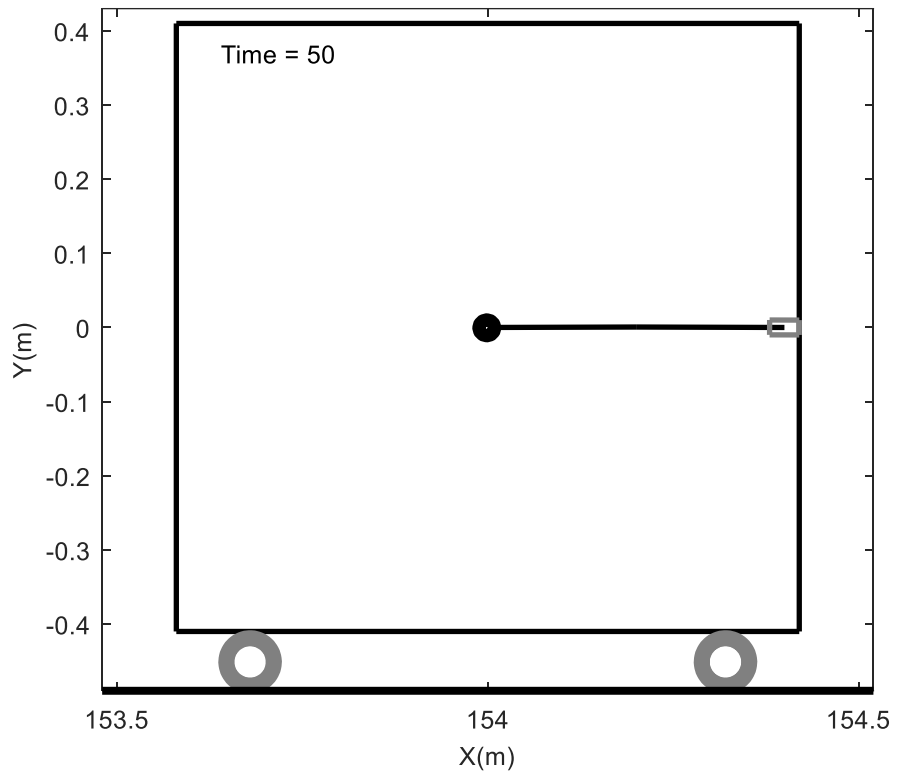

(e)

**Fig. 6** Snapshots of the system's motion of the first case study

*3.2. An unconstrained multibody system including three connected rigid bodies translating and rotating in space*

As the second case study, the system shown in Figure 7 is analyzed. This system is a set of three connected rigid bodies. In addition to three external torques applied on the main body (the middle body), two motors are mounted on joints $A$ and $B$ to produce $\boldsymbol{\tau}_{p1}$ and $\boldsymbol{\tau}_{p2}$ respectively. The purpose of this example is to show the application of the proposed method in unconstrained problems.

The GCs vector is chosen as:

$$\boldsymbol{q} = [\psi \quad \theta \quad \varphi \quad \gamma_1 \quad \gamma_2 \quad X \quad Y \quad Z]^T \tag{72}$$

in which $\psi$, $\theta$ and $\varphi$ are the amounts of the sequence of rotations that $xyz$ frame undergoes about its own axes ($z$, $y$ and $x$, respectively) to transform from its initial orientation to its final orientation. Furthermore, $X$, $Y$ and $Z$ determine the absolute position of the center of mass of the middle body. Other parameters are shown in Figure 7. The system has no constraints, so it has $p = m = 8$ number of DOFs. After deriving generalized forces (similar to Eq. (70)) it turns out that based on *Definition1*, the ignorable coordinates vector is:

$$\boldsymbol{q}_I = [X \quad Y \quad Z]^T \tag{73}$$

Hence, the system has $s = 3$ number of ignorable coordinates. The vector of quasi-velocities is selected as:

$$\boldsymbol{u} = [\boldsymbol{\omega}_m{}^T \quad \dot{\gamma}_1 \quad \dot{\gamma}_2 \quad \dot{X} \quad \dot{Y} \quad \dot{Z}]^T \tag{74}$$

Where $\boldsymbol{\omega}_m$ is the angular velocity of the middle body relative to the $xyz$ frame. The first five elements of $\boldsymbol{u}$ construct the $\boldsymbol{u}_{NI}$ vector which is to be used in the proposed method.

Similar to the first example, all external torques are set to zero so that the mechanical energy of the system remains constant. The inertia matrix of the middle body is:

$$\boldsymbol{I} = \begin{bmatrix} I_{xx} & -I_{xy} & -I_{xz} \\ -I_{xy} & I_{yy} & -I_{yz} \\ -I_{xz} & -I_{yz} & I_{zz} \end{bmatrix} \tag{75}$$

For the sake of simplicity in calculating inertia properties of other two bodies, they are considered as uniform rectangular planes.

Initial conditions and properties of the bodies of the system are provided in Tables 4 and 5.

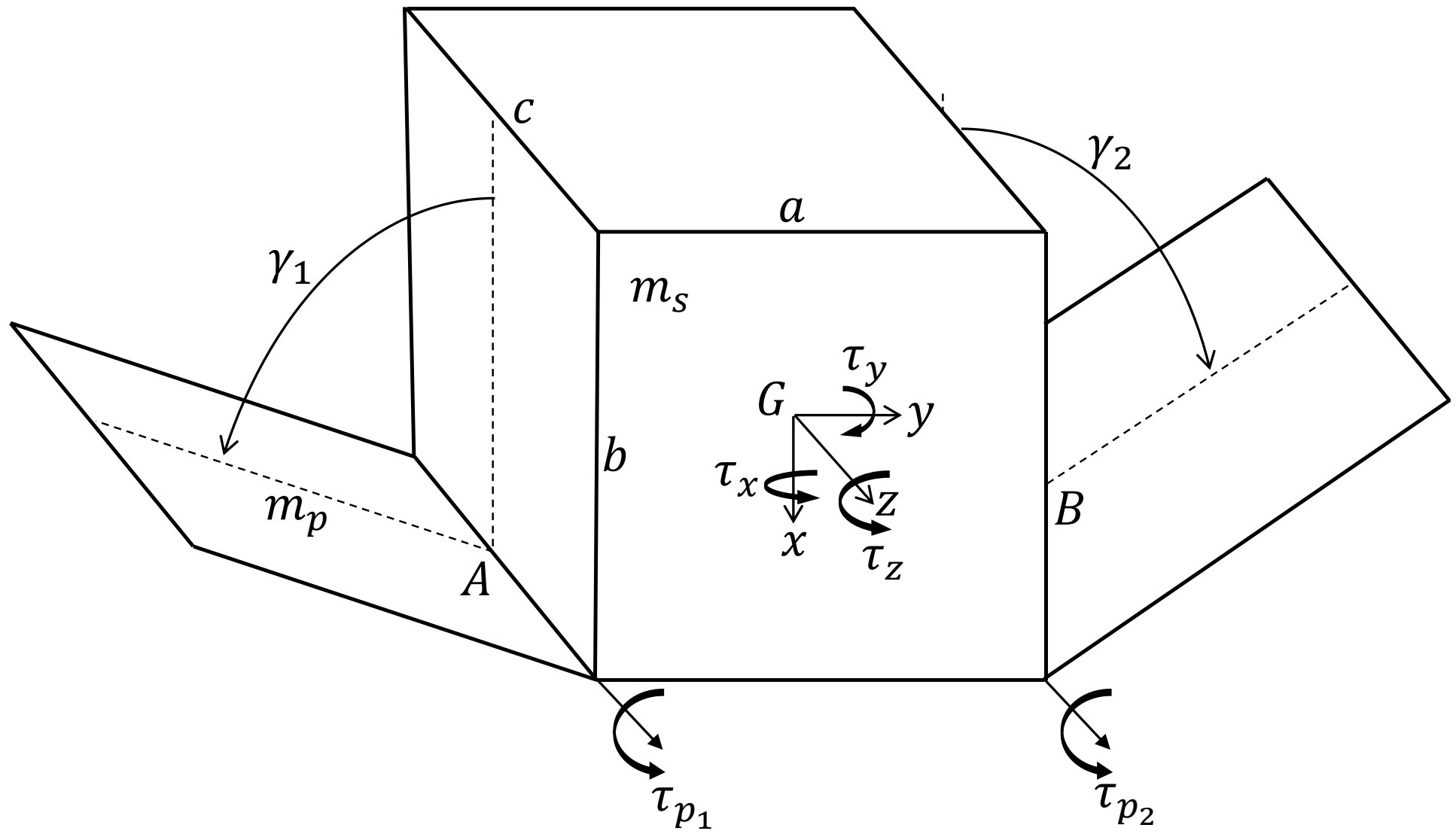

**Fig. 7** Second simulation example: a system of three connected rigid bodies

**Table 4.** Parameters of the system considered in the second example [41]

| ***Parameter*** | ***Value*** | ***Unit*** |
|---|---|---|
| $m_s$ | 100 | $kg$ |
| $m_p$ | 10 | $kg$ |
| $a$ | 2 | $m$ |
| $b$ | 2 | $m$ |
| $c$ | 2 | $m$ |
| $I_{xx}$ | 67 | $kgm^2$ |
| $I_{yy}$ | 67 | $kgm^2$ |
| $I_{zz}$ | 67 | $kgm^2$ |
| $I_{xy}$ | 5 | $kgm^2$ |
| $I_{xz}$ | 2 | $kgm^2$ |
| $I_{yz}$ | 0 | $kgm^2$ |

**Table 5.** Initial conditions of the system considered in the second example [41]

| ***Parameter*** | ***Value*** | ***Unit*** |
|---|---|---|
| $\psi$ | $\pi/10$ | $rad$ |
| $\theta$ | $\pi/6$ | $rad$ |
| $\varphi$ | 0.03 | $rad$ |
| $\gamma_1$ | 0 | $rad$ |
| $\gamma_2$ | 0 | $rad$ |
| $X$ | 3 | $m$ |
| $Y$ | 3 | $m$ |
| $Z$ | 9 | $m$ |
| $\dot{\psi}$ | 0.3 | $rad/s$ |
| $\dot{\theta}$ | 0 | $rad/s$ |
| $\dot{\varphi}$ | -0.3 | $rad/s$ |
| $\dot{\gamma}_1$ | 0 | $rad/s$ |
| $\dot{\gamma}_2$ | 0.2 | $rad/s$ |
| $\dot{X}$ | 0.1 | $m/s$ |
| $\dot{Y}$ | -0.3 | $m/s$ |
| $\dot{Z}$ | -0.4 | $m/s$ |

Similar to what was done in the first example, the problem is solved via 4 different methods: the Lagrange method, the Maggi method, the Kane (or Gibbs-Appell) method, and the proposed method. The simulation is conducted for 50 seconds, the time step is fixed at 0.1 (s), and similar to the prior example, the ode45 function of MATLAB is employed to solve the equations of motion. The results are provided in Table 6. It's worthy to note that although linear momentum conservation error was computed in $X$, $Y$ and $Z$ directions, Table 6 only represents this error in $X$ direction. This is because the order of the error in all three directions was the same for each approach.

As it can be seen from Table 6, the proposed method acted significantly better than other approaches. For the ease of comparison, graphs of energy and linear momentum (in $X$ direction) errors are provided in Figures 8 and 9. Again, it can be easily comprehended that the proposed method conserves mechanical energy and linear momentum noticeably better than other methods. The plot of the generalized coordinates is shown in Figure 10. Furthermore, in order to exhibit the motion, snapshots of the motion of the system in several instants are displayed in Figure 11.

**Table 6.** Results of simulation of the second example using different methods

| ***Method*** | Lagrange | Maggi | Kane (Gibbs-Appell) | Proposed |
|---|---|---|---|---|
| ***Number of State Variables*** | 16 | 16 | 16 | 13 |
| ***Number of Equations*** | 8 | 8 | 8 | 5 |
| ***CPU Time (s)*** | 1.1 | 0.99 | 0.68 | 0.67 |
| ***Norm of Energy Error*** | $2.96 \times 10^{-7}$ | $2.96 \times 10^{-7}$ | $9.10 \times 10^{-9}$ | $6.36 \times 10^{-9}$ |
| ***Norm of Linear Momentum Error (X)*** | $1.55 \times 10^{-6}$ | $1.55 \times 10^{-6}$ | $3.75 \times 10^{-8}$ | 0 |

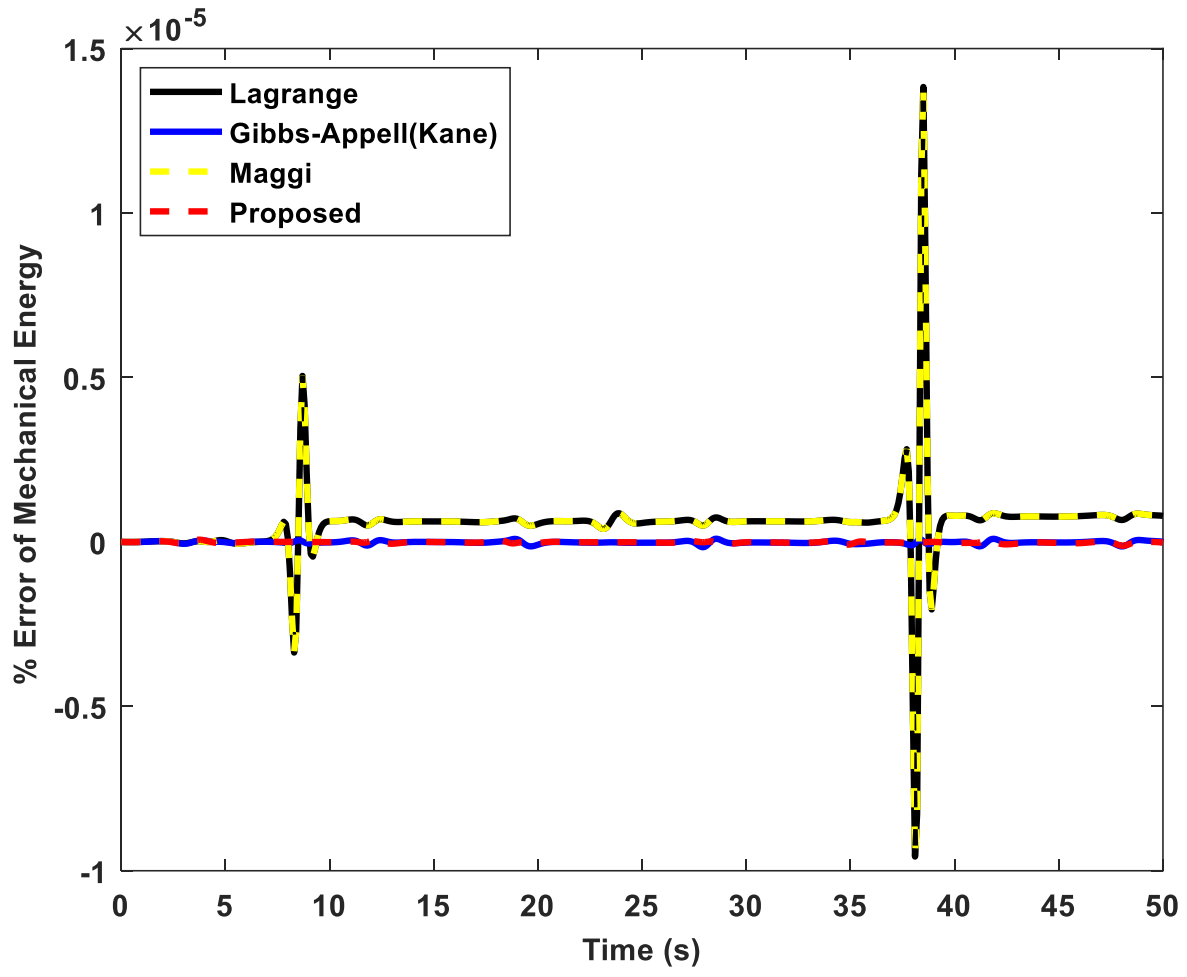

**Fig. 8** Percentage of mechanical energy drift in the second case study

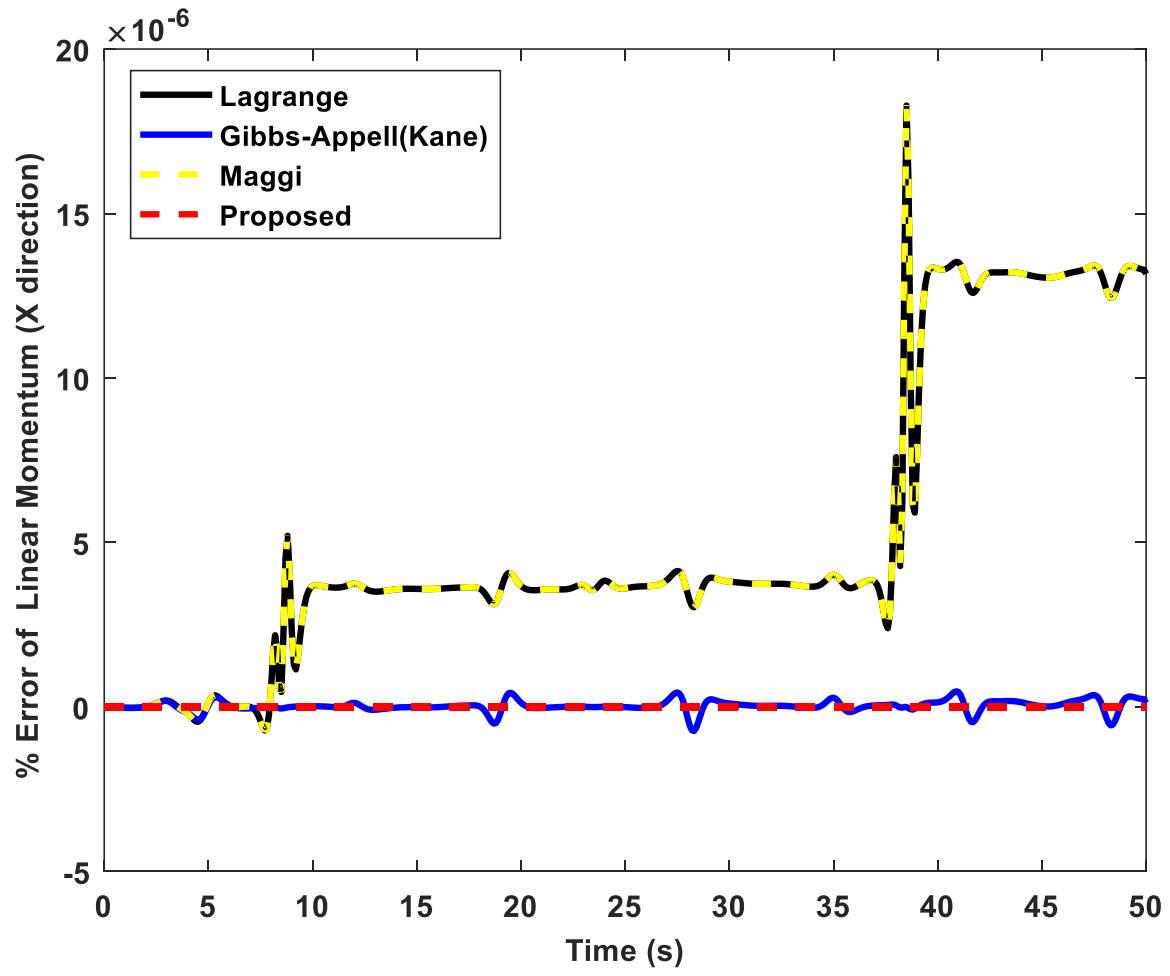

**Fig. 9** Percentage of linear momentum (x direction) conservation error in the second case study

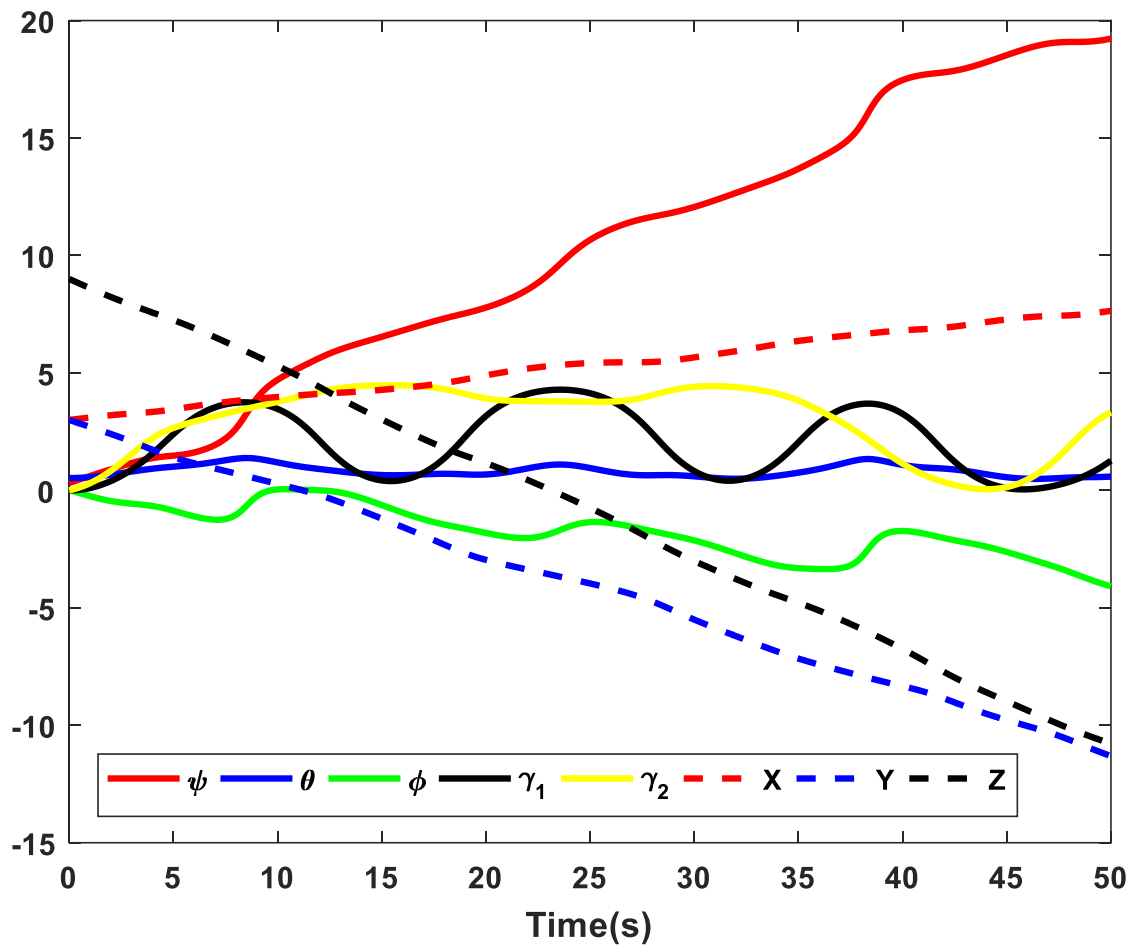

**Fig. 10** The plot of the generalized coordinates in the second case study

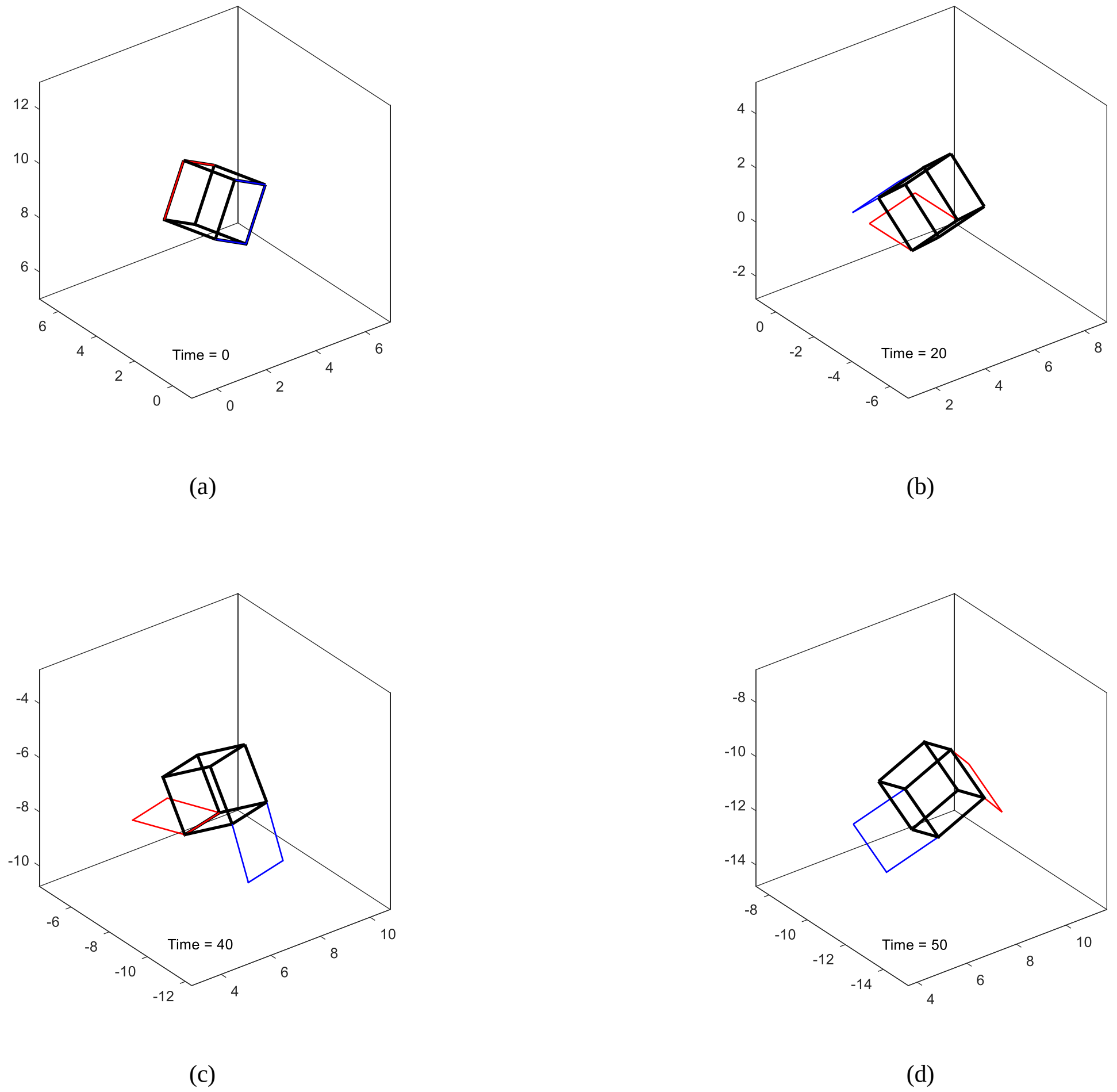

**Fig. 11** Snapshots of the system's motion of the second case study

### *3.3. A cubic satellite on which a deployable boom is mounted*

As the last case study, a cubic satellite is considered. On one side of this satellite a deployable boom is mounted which expands and thus moves a mass to a predetermined position. The system is shown in Figure 12. A motor produces force $\boldsymbol{F}$, thus pushing the mass $m$ forward. The initial length of the boom is $a$. The mass of the boom which is used to move the mass $m$ is negligible.

The GCs vector is chosen as:

$$\boldsymbol{q} = [\psi \quad \theta \quad \varphi \quad \rho \quad X \quad Y \quad Z]^T \tag{76}$$

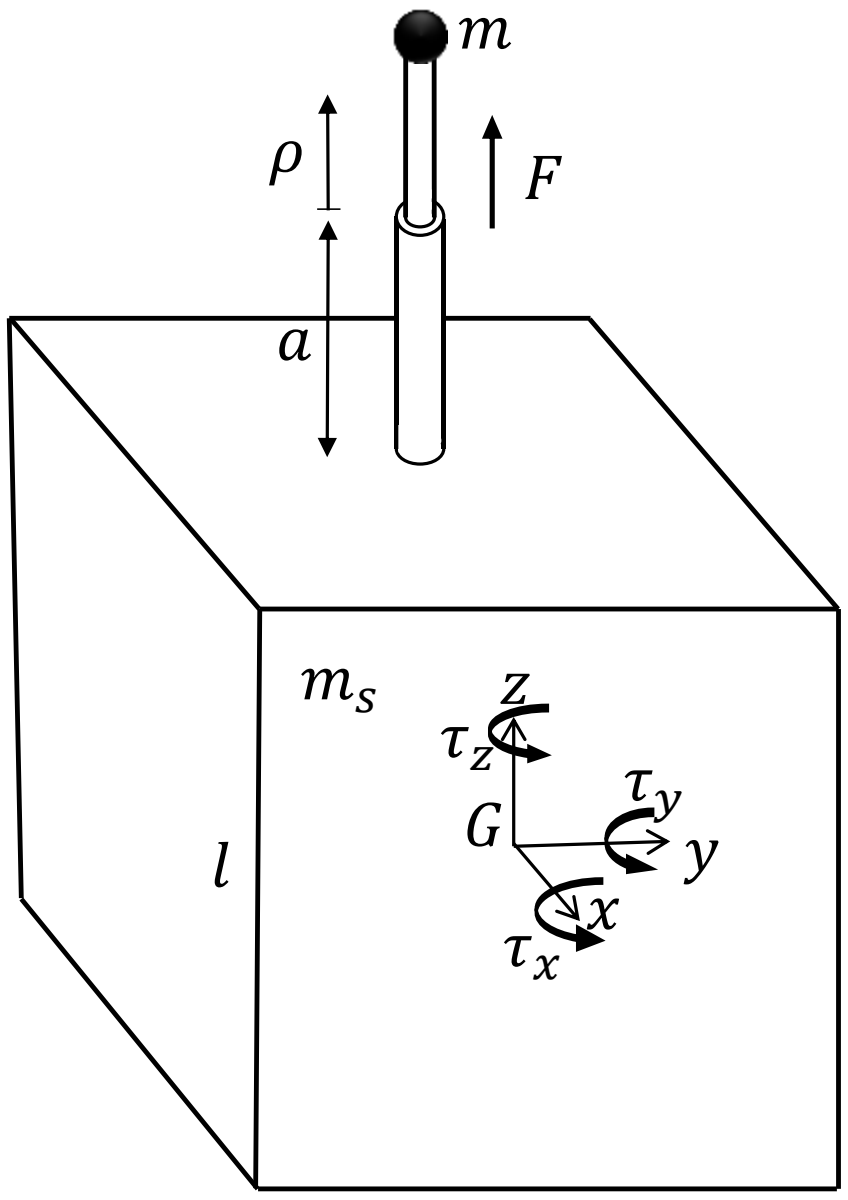

**Fig. 12** Third simulation example: a satellite with deployable magnetometer boom

The first and the last three parameters used in Eq. (76) are defined as same as those used in Eq. (72). $\rho$ is the amount of variation in the length of the boom.

The system has no constraints, so it has $p = m = 7$ number of DOFs. Similar to the previous example, after deriving generalized forces (similar to Eq. (70)) it turns out that based on *Definition1*, the vector of ignorable coordinates is:

$$\boldsymbol{q}_I = [X \quad Y \quad Z]^T \tag{77}$$

Hence, the system has $s = 3$ number of ignorable coordinates. The vector of quasi-velocities is selected as:

$$\boldsymbol{u} = [\boldsymbol{\omega}_m{}^T \quad \dot{\rho} \quad \dot{X} \quad \dot{Y} \quad \dot{Z}]^T \tag{78}$$

where $\boldsymbol{\omega}_m$ is the angular velocity of the satellite relative to the $xyz$ frame. The first four elements of $\boldsymbol{u}$ construct the $\boldsymbol{u}_{NI}$ vector which is to be used in the proposed method.

Similar to the previous example, all external torques are set to zero. However, the magnitude of force $\boldsymbol{F}$ changes over time as:

$$|\boldsymbol{F}| = -0.018\, sin(0.089t) + 0.012\, cos(0.0485t) \tag{79}$$

Accordingly, unlike previous examples, mechanical energy does not remain constant. Instead, the amount of energy error should remain constant. The energy error is defined as:

$$Energy\,Error = T + V - \int_0^{t_s} \boldsymbol{Q}^T \dot{\boldsymbol{q}}\, dt \tag{80}$$

in which $T$ and $V$ are the kinetic and potential energies of the system, $\dot{\boldsymbol{q}} \in \mathbb{R}^m$ is the generalized velocity vector, $\boldsymbol{Q} \in \mathbb{R}^m$ is the generalized force vector expressed in the generalized velocity space, and $t_s$ is the simulation time. Moreover, the term $\boldsymbol{Q}^T \dot{\boldsymbol{q}}$ indicates the input power to the system. For approaches that employ the quasi-velocity concept, the equivalent term for the power will be $\boldsymbol{U}^T \boldsymbol{u}$, in which $\boldsymbol{u}$ is the quasi-velocity vector, and $\boldsymbol{U}$ is the generalized force vector expressed in the quasi-velocity space. The inertia matrix of the satellite is the same as Eq. (75).

The system parameters and initial conditions and are listed in Tables 7 and 8, respectively.

**Table 7.** Parameters of the system considered in the third example

| ***Parameter*** | ***Value*** | ***Unit*** |
|---|---|---|
| $m_s$ | 2000 | $kg$ |
| $m$ | 1 | $kg$ |
| $l$ | 2.5 | $m$ |
| $a$ | 0.5 | $m$ |
| $I_{xx}$ | 1400 | $kgm^2$ |
| $I_{yy}$ | 900 | $kgm^2$ |
| $I_{zz}$ | 1100 | $kgm^2$ |
| $I_{xy}$ | 5 | $kgm^2$ |
| $I_{xz}$ | 8 | $kgm^2$ |
| $I_{yz}$ | 3 | $kgm^2$ |

**Table 8.** Initial conditions of the system considered in the third example

| ***Parameter*** | ***Value*** | ***Unit*** |
|---|---|---|
| $\psi$ | $\pi/8$ | $rad$ |
| $\theta$ | $\pi/12$ | $rad$ |
| $\varphi$ | 0.08 | $rad$ |
| $\rho(= a)$ | 0.5 | $m$ |
| $X$ | 0 | $m$ |
| $Y$ | 0 | $m$ |
| $Z$ | 0 | $m$ |
| $\dot{\psi}$ | -0.1 | $rad/s$ |
| $\dot{\theta}$ | 0.05 | $rad/s$ |
| $\dot{\varphi}$ | -0.05 | $rad/s$ |
| $\dot{\rho}$ | 0 | $m/s$ |
| $\dot{X}$ | 2 | $m/s$ |
| $\dot{Y}$ | 1 | $m/s$ |
| $\dot{Z}$ | 0 | $m/s$ |

Similar to what was done in the last two examples, the problem is solved via 4 different methods: Lagrange method, Maggi method, Kane (or Gibbs-Appell) method, and the proposed method. The simulation is conducted for 50 seconds, the time step is fixed at 0.1 (s), and similar to the previous examples, the ode45 function of MATLAB is employed to solve the equations of motion. The results are provided in Table 9. It's worthy to note that although linear momentum conservation error was computed in $X$, $Y$ and $Z$ directions, Table 9 only represents this error in $X$ direction. This is because the order of the error in all three directions was the same for each approach.

As it can be seen from the Table 9, the proposed method acted significantly better than other approaches. For the ease of comparison, graphs of power and linear momentum (in $X$ direction) errors are provided in Figures 13 and 14. Again, it can be easily understood that the proposed method conserves power and linear momentum noticeably better than other methods. The plot of the generalized coordinates is shown in Figure 15. Similar to previous examples, shots of the motion of the system in several random instants are exhibited in Figure 16.

**Table 9.** Results of simulation of the third example employing various methods

| ***Method*** | Lagrange | Maggi | Kane (Gibbs-Appell) | Proposed |
|---|---|---|---|---|
| ***Number of State Variables*** | 15 | 15 | 15 | 12 |
| ***Number of Equations*** | 7 | 7 | 7 | 4 |
| ***CPU Time (s)*** | 0.59 | 0.53 | 0.50 | 0.48 |
| ***Norm of Energy Error*** | $6.31 \times 10^{-14}$ | $6.18 \times 10^{-14}$ | $5.01 \times 10^{-14}$ | $7.91 \times 10^{-15}$ |
| ***Norm of Linear Momentum Error (X)*** | $3.23 \times 10^{-14}$ | $3.75 \times 10^{-14}$ | $3.74 \times 10^{-14}$ | 0 |

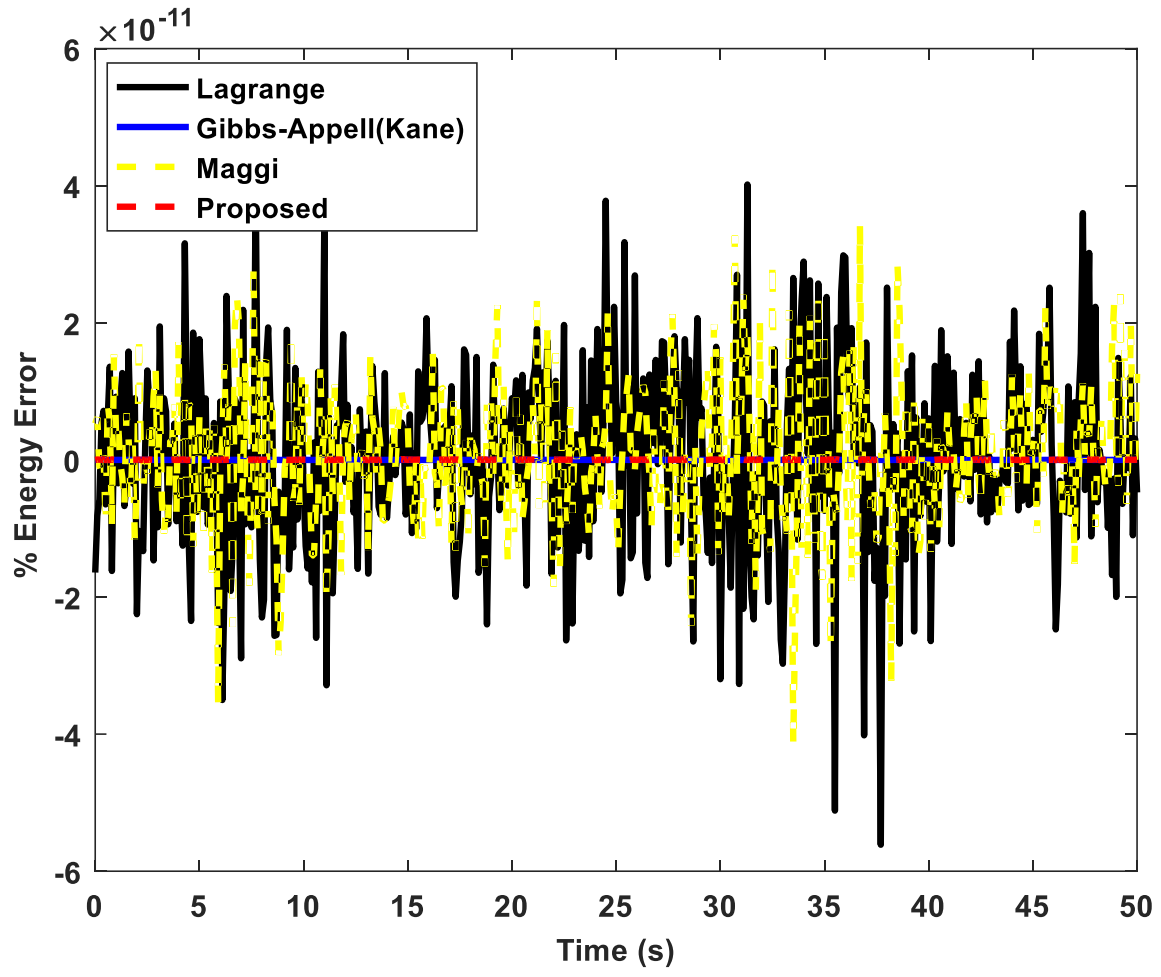

**Fig. 13** Percentage of energy error in the third case study

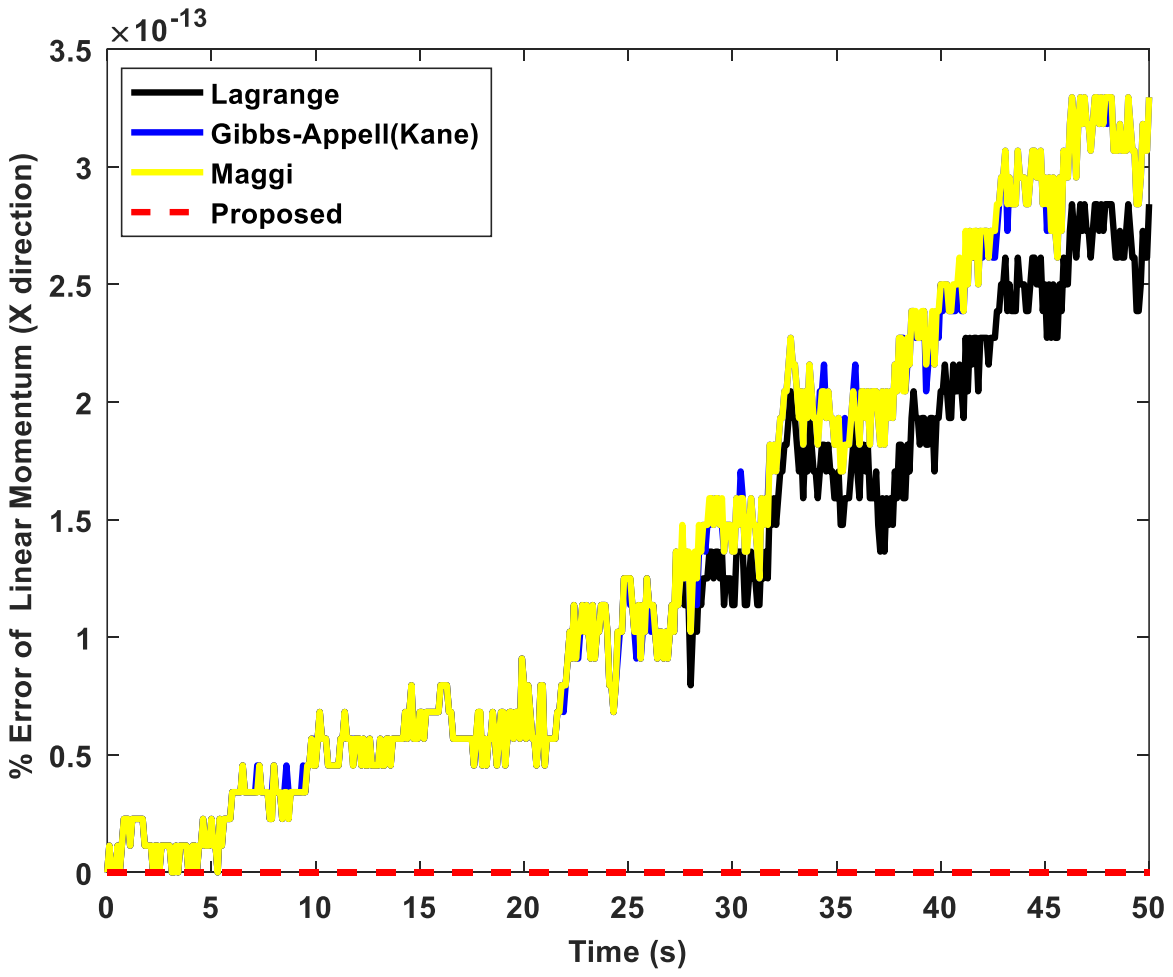

**Fig. 14** Percentage of linear momentum (x direction) conservation error in the third case study

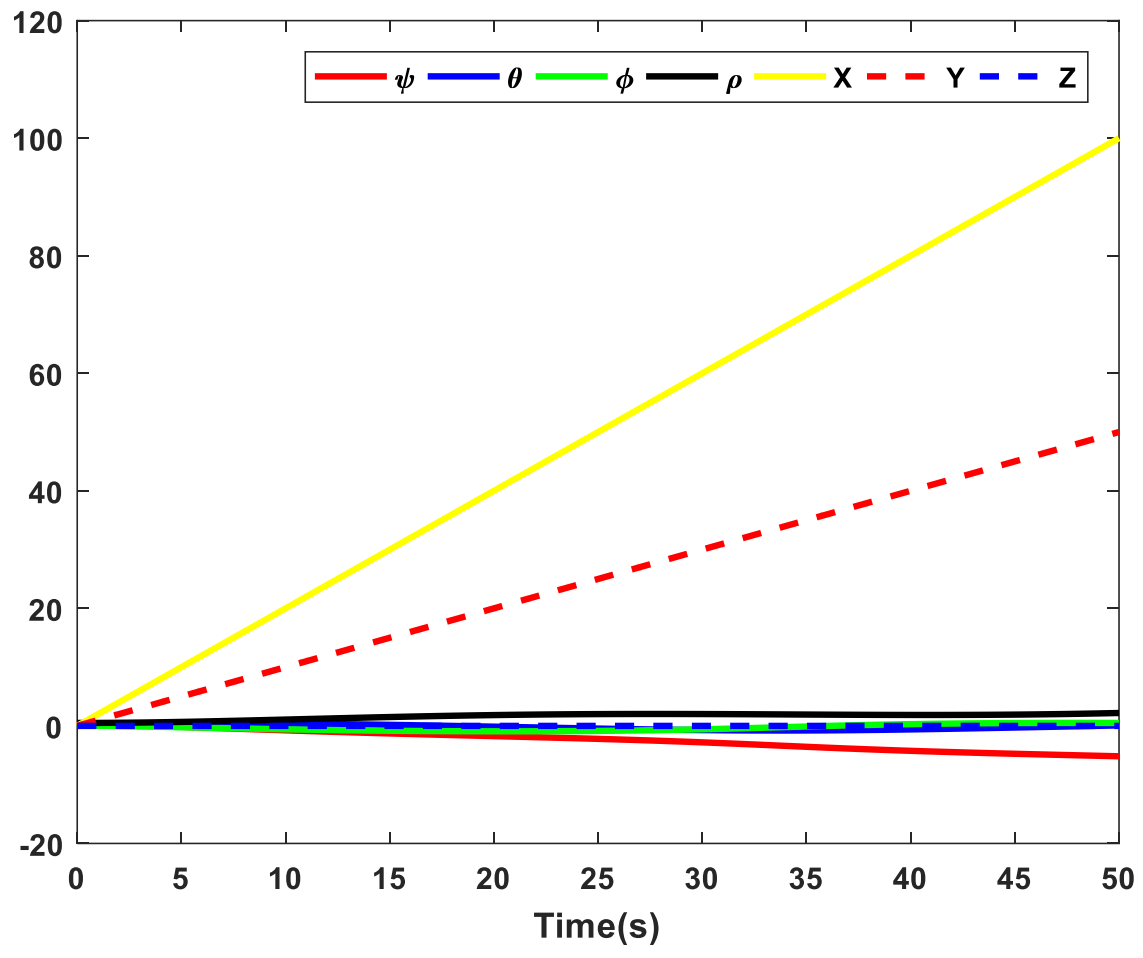

**Fig. 15** The plot of the generalized coordinates in the third case study

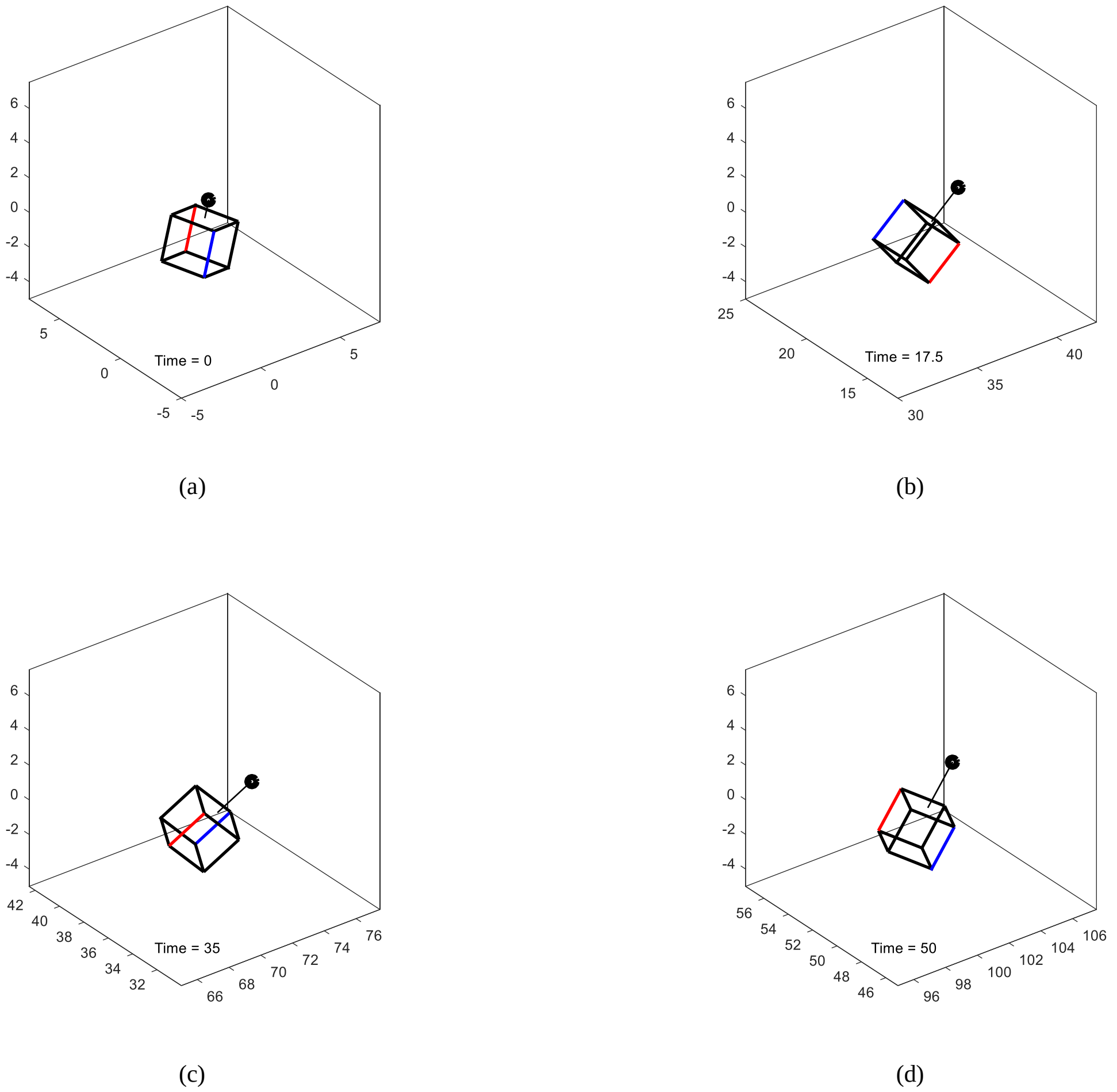

**Fig. 16** Snapshot of motion of the system analyzed as the third case study

## 4. Conclusion

In this paper, we presented an efficient form of Volterra's equations of motion for systems including ignorable coordinates. Meanwhile, we proposed the concept of "dynamical constraints" associated with the definition of ignorable coordinates. Subsequently, by augmenting dynamical constraints along with ordinary constraints (which were regarded as "kinematical constraints") with the equations derived from the proposed method, the minimum number of equations required for complete motion analysis were derived. Three simulation examples were provided to show the usage of the approach. The first one was a constrained system including a cart and a 2-DOF pendulum. The problem was solved using 4 different methods, that is, Maggi, Gibbs-Appell (or Kane), Lagrange, and the proposed method. In the second example, an unconstrained multibody system including three connected rigid bodies with complex motion in space was considered. The third problem comprised a cubic satellite on which a deployable boom is attached. This boom is used to move a mass to a desired position. These problems were solved using the same methods that were employed in the first example. In all these three case studies, the better performance of the proposed method, in terms of CPU time and computational error, was demonstrated. The developed method can be generalized and used in practical large-scale systems including ignorable coordinates and result in significant reduction in computation time and violation of motion constants, when compared to conventional methods. In conclusion, the contributions of this paper are:

- Introducing dynamical constraints and embedding them with motion equations
- Developing a novel Volterra-based approach to derive equations of motion, which outperforms existing methods in terms of runtime and constraints violation
- Possessing the minimum number of equations for systems including ignorable coordinates

## 5. Nomenclature

| | |
|---|---|
| $\boldsymbol{a}$ | Constraints' Jacobian matrix |
| $\boldsymbol{b}$ | Constraints' bias vector |
| $\boldsymbol{q}_{NI}$ | Non-ignorable generalized coordinates vector |
| $\boldsymbol{q}_{I}$ | Ignorable generalized coordinates vector |
| $\boldsymbol{H}$ | Angular momentum vector |
| $\boldsymbol{P}$ | Linear momentum vector |
| $\boldsymbol{\omega}$ | Angular velocity |
| $\boldsymbol{V}$ | Linear velocity |
| $\boldsymbol{u}$ | Quasi-velocity vector |
| $\boldsymbol{Q}$ | Generalized forces vector |
| $\boldsymbol{U}$ | Generalized forces vector in quasi-velocity space |
| $\boldsymbol{G}$ | Generalized momentum vector |
| $\boldsymbol{F}$ | External force vector |
| $\boldsymbol{Y}$ | Jacobian matrix of quasi-velocities |
| $\boldsymbol{Z}$ | Bias vector of quasi-velocities |
| $T$ | Kinetic energy |
| $\boldsymbol{M}$ | Mass matrix |
| $\boldsymbol{I}$ | Inertia matrix |
| $\mathcal{L}$ | Lagrangian function |